\documentclass{article}
\usepackage{authblk}

\usepackage{geometry}
\usepackage[title]{appendix}
\usepackage{graphicx}
\usepackage{mathptmx}      
\usepackage{multirow}
\usepackage{tikz}
\usetikzlibrary{positioning}
\usepackage[hidelinks]{hyperref}
\usepackage{booktabs, array}
\usepackage{eqnarray, amsmath}
\usepackage{xcolor}
\usepackage[labelfont=bf,labelsep=space,font={small}]{caption}
\usepackage{subcaption}
\usepackage{tabu}
\usepackage[scr=rsfso,cal=zapfc,frak=euler,bb=ams]{mathalfa}
\usepackage{natbib}
\newcommand{\mike}[1]{\textcolor{black}{#1}}

\setcitestyle{round, aysep={ },yysep={;}}

\title{NASH EQUILIBRIA IN CERTAIN TWO-CHOICE MULTI-PLAYER GAMES PLAYED ON THE LADDER GRAPH}
\author{}

\makeatletter
\renewcommand\@date{{%
  \vspace{-\baselineskip}%
  \large\centering
  \begin{tabular}{@{}c@{}}
    Victoria Sánchez Muñoz\textsuperscript{1} \\
    \normalsize v.sanchezmunoz1@nuigalway.ie
  \end{tabular}%
  \quad and\quad
  \begin{tabular}{@{}c@{}}
    Michael Mc Gettrick\textsuperscript{2} \\
    \normalsize michael.mcgettrick@nuigalway.ie
  \end{tabular}

  \bigskip
  \textsuperscript{1,2}School of Mathematics, Statistics and Applied Mathematics\\ National University of Ireland Galway\\
Galway, Ireland\par
  \bigskip
  \bigskip
}}
\makeatother

\begin{document}

\maketitle

\begin{abstract}
In this article we \mike{compute analytically} the number of Nash Equilibria (NE) for a two-choice game played on a \mike{(circular) ladder graph} with $2n$ players. We \textcolor{black}{consider a set of games with generic payoff parameters, with the only} requirement that a NE occurs if the players choose opposite strategies (anti-coordination game). The results show that for both, the ladder and circular ladder, the number of NE grows exponentially with (half) the number of players $n$, as $N_{NE}(2n)\sim C(\varphi)^n$, where $\varphi=1.618..$ is the golden ratio and $C_{circ}>C_{ladder}$. In addition, the value of the scaling factor $C_{ladder}$ depends on the value of the payoff parameters. However, that is no longer true for the circular ladder (3-degree graph), that is $C_{circ}$ is constant, which might suggest that the topology of the graph indeed plays an important role for setting the number of NE.
\end{abstract}

\bigbreak
{\emph{Keywords}: graphical game; ladder graph; circular ladder; Nash Equilibrium}
\bigbreak
{Subject Classification: 91A43, 05C57.}

\section{Introduction}
\label{sec:intro}
Graphical games [\citet{Kearns07chapter7}] is an area within game theory with several applications in biology (with evolutionary graphical games [\citet{gamesongraphs,Allen2019}]), economics [\citet{Leduc_2017,strategicinvesvestmentgraphgames}], sociology [\citet{Eger_2016}] and computer science [\citet{Cibulka_2013,Gottlob_2005}]. A graph $\mathcal{G}$ and a predefined game is considered, in which the vertices of the graph $\mathcal{G}$ represent the players and each player plays a game with their neighbours, that is, the players connected to them. 

\textcolor{black}{The main focus on graphical games lies heavily on computational game theory, i.e. how to find (pure and/or mixed) NE solutions [\citet{Nash48}] from an algorithmic point of view (see for example [\citet{graphmodelsforGT}] in which a tree graph is considered) and the complexity class of those, as done for instance in [\citet{PureNEinGraphGameCSA}], [\citet{complexityCNE}] and [\citet{Dyer2011PairwiseInteractionG}].
Other authors address more game-theoretical questions as to whether any deterministic graphical multi-person game has a NE in pure stationary strategies. To address the question they may give and analyse ``easy" examples of graphs with a low number of players to disprove conjectures that would answer in the positive (see for example [\citet{BOROS201821}]). However, throughout our paper, our approach is completely different. We compute analytically an expression to get the number of all pure NE of $2n$ players who pair-wise play a generic anti-coordination game, whose connection to the opponents follows two particular and well-known graphs: the ladder and the circular ladder. }

\textcolor{black}{
\mike{One of our motivations is to
extend results for
graphical games on 2-regular connected graphs
(i.e. cycles) to a family of 3-regular connected
graphs.} Our way of proceeding might be a basis when considering other (highly regular and possibly complete) graphs, as well as \mike{analytically} presenting results to check with the algorithms.}
\mike{ While we study the asymptotic growth of NE as the number of players increases (and hence the graph increases), other authors fix the number of players, and study the asymptotic growth of NE as the number of strategy choices increases (see [\citet{MCLENNAN2005264}]). It is very interesting to observe that in both cases, exponential growth is obtained. Of  course, both models overlap and can be compared for 2-player 2-choice games.}

\section{Definition of the game}\label{sec:defgame}
Our underlying game is a non-cooperative, non-zero sum, two-choice game. Considering two players: Player 1's strategy set is denoted by $S_1=\{a,b\}$, and a specific strategy $s_1 \in S_1$, and similarly for Player 2, $S_2=\{a,b\}$, with $s_2 \in S_2$.

The typical payoff matrix is in table \ref{tab:payoffmatrix}, whose four real parameters $p$, $r$, $q$ and $s$ characterise which (if any) NE solution/s there are in this two-player scenario. Usually, two conditions are imposed with such parameters so that the NE strategies are $\{s_1=a, s_2=b\}$ and its symmetric counterpart $\{s_1=b, s_2=a\}$. These two conditions are: $q>s$ and $r>p$. Examples of such games include \emph{the game of chicken} or any anti-coordination game, in which the best response is to play the opposed action to your opponent's.

\begin{table}
	\caption{The first value in brackets denotes Player 1's payoff, and the second Player 2's. The bold values correspond to the payoffs of the typical NE strategies.}
    \centering
    \begin{tabular}{cc|cc|} 
        \cline{3-4}
        & & \multicolumn{2}{c|} {Player 2} \\
        &  &  a  &  b \\ \cline{1-4}
     \multicolumn{1}{|c}{} \multirow{2}{*}{Player 1}  & a &  (p, p)  & (\textbf{q}, \textbf{r})  \\
      \multicolumn{1}{|c}{} & b & (\textbf{r}, \textbf{q}) & (s, s) \\
      \hline
    \end{tabular}
    
    \label{tab:payoffmatrix}
\end{table}

\begin{figure}
\centering
\begin{tikzpicture}
\draw[black, thick] (0,0) -- (3.7,0);
\draw[black, thick] (4.3,0) -- (8,0);
\draw[black, thick] (0,1) -- (3.7,1);
\draw[black, thick] (4.3,1) -- (8,1);
\draw[black, thick] (0,0) -- (0,1) ;
\draw[black, thick] (1,0) -- (1,1) ;
\draw[black, thick] (2,0) -- (2,1) ;
\draw[black, thick] (3,0) -- (3,1) ;
\draw[black, thick] (5,0) -- (5,1) ;
\draw[black, thick] (6,0) -- (6,1) ;
\draw[black, thick] (7,0) -- (7,1) ;
\draw[black, thick] (8,0) -- (8,1) ;
\filldraw[color=black, thin] (3.8,0.5) circle (0.65pt);
\filldraw[color=black, thin] (4.0,0.5) circle (0.65pt);
\filldraw[color=black, thin] (4.2,0.5) circle (0.65pt);
\filldraw[color=black] (0,0) circle (2pt);
\filldraw[color=black] (1,0) circle (2pt);
\filldraw[color=black] (2,0) circle (2pt);
\filldraw[color=black] (3,0) circle (2pt);
\filldraw[color=black] (5,0) circle (2pt);
\filldraw[color=black] (6,0) circle (2pt);
\filldraw[color=black] (7,0) circle (2pt);
\filldraw[color=black] (8,0) circle (2pt);
\filldraw[color=black] (0,1) circle (2pt);
\filldraw[color=black] (1,1) circle (2pt);
\filldraw[color=black] (2,1) circle (2pt);
\filldraw[color=black] (3,1) circle (2pt);
\filldraw[color=black] (5,1) circle (2pt);
\filldraw[color=black] (6,1) circle (2pt);
\filldraw[color=black] (7,1) circle (2pt);
\filldraw[color=black] (8,1) circle (2pt);
\end{tikzpicture}
\caption{Ladder graph. Each node represents a player, who will play a game with the players connected to them. The ending players only play two games, while the others play three games.}
\label{fig:openladd}
\end{figure}
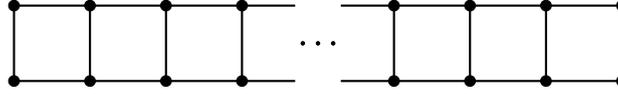

Now that we have defined that two-player game, we can extend this to a game played in the ladder graph, in which each player, represented as a vertex, plays a game with all the players that are connected to them. We will denote the total number of players as $2n$. Figure \ref{fig:openladd} shows the set-up.

Following the convention, let's denote the strategy set of Player $i$ by $S_i=\{a,b\}$, $s_i$ as Player $i$'s pure strategy ($s_i \in S_i$), $S_{-i}$ as the strategy set of Player $i$ 's opponents, $s_{-i}$ as an element of  $S_{-i}$, and $s^*_i$ to the best response of Player $i$ given $s_{-i}$. 

In this ladder set-up, there are two distinct type of players, the ones at both ends of the ladder and the middle players. The ending players only play two games, while the middle ones play three games. Let us focus on the first player Player 1, located at the start of the top row of the ladder, their strategy set is $S_1=\{a,b\}$, and the strategy set of Player 1's opponents is $S_{-1}=\{aa, ab, ba, bb\}$, in which $xy$, with $x,y\in\{a,b\}$, means that the player to the right of Player 1 chooses $x$, while the one below chooses $y$. Nevertheless, if we look at some player in the middle, for example Player 3 (to the right of Player 1), their strategy set remains $S_3=\{a,b\}$ and their opponent's strategy set is $S_{-3}=\{aaa, aab, aba, baa, abb, bab, bba, bbb\}$, in which $xyz$, with $x,y,z\in\{a,b\}$, means that the player to the left of Player 3 chooses $x$, while the one below chooses $y$ and the one to the right of Player 3 plays $z$.

We want this game to be fair, so to equalise every player, the payoffs for each player will be averaged. Table \ref{tab:players} shows all the payoffs for Player $i$ when choosing to play $s_i=a$ or $s_i=b$ against all the different\footnote{We ignore the elements in $S_{-i}$ that are just a permutation of another element because the payoff for Player $i$ would remain unchanged.} possible strategies of their opponents.

\begin{table}[]
    \centering
    \caption{Payoff values for Player $i$. The upper half of the table corresponds to a player who is in the middle of the ladder and therefore, plays three games, while the lower half represents a player in the ending of the ladder (plays only two games).}
    \begin{tabular}{cc|cc|}
       \cline{3-4} & &  \multicolumn{2}{c|}{Player $i$'s strategy $s_i$} \\
         & & a & b \\ \hline
       \multicolumn{1}{|c}{} \multirow{7}{*}{$s_{-i}$} & aaa  & $p$ & $\textbf{r}$ \\
       \multicolumn{1}{|c}{} & bbb  & $\textbf{q}$ &  $s$ \\
       \multicolumn{1}{|c}{} & aab  & $(2p+q)/3$ & $(2r+s)/3$ \\
       \multicolumn{1}{|c}{} & abb  & $(2q+p)/3$ & $(2s+r)/3$ \\ \cline{3-4}
       \multicolumn{1}{|c}{} & aa & $p$ & $\textbf{r}$ \\
       \multicolumn{1}{|c}{} & bb & $\textbf{q}$ &  $s$\\
       \multicolumn{1}{|c}{} & ab & $(p+q)/2$ & $(r+s)/2$\\ \hline
    \end{tabular}
    \label{tab:players}
\end{table}

The bold values are, by the assumption of the two-player game, greater than the other values in the same row. Hence, that would set the preferred strategy for Player $i$ given $s_{-i}$. For instance, if $s_{-i}=aaa$, then $s^{*}_{i}=b$. Nevertheless, for the case when $s_{-i}=ab$ (the player is at one end of the ladder), $s_{-i}=aab$ and $s_{-i}=abb$ for the middle player, the highest value of the row is undefined. That means, we need to make an additional assumption in order to define the best strategy for Player $i$ in those situations.
To do that, we define a function as the difference in the payoffs between the case in which $s_i=a$ and the case in which $s_i=b$ given $s_{-i}$. The study of the sign of such function sets the conditions in which Player $i$ would rather play $a$ ($b$) than playing $b$ ($a$) against that particular $s_{-i}$. There are three different undefined situations, so we need three functions. The first two functions are meant for the middle players and the last one, for the players at the beginning/ending of the ladder. The three functions are\footnote{Since we are only interested in the sign of the function, we ignore the common denominators in the payoffs.}:

\begin{eqnarray}\label{fgh0}
    & f(p,q,r,s)=(2p+q)-(2r+s)=-2(r-p)+(q-s) \\
    & g(p,q,r,s)=(2q+p)-(2s+r)=2(q-s)-(r-s) \\
    & h(p,q,r,s)=(p+q)-(r+s)=(q-s)-(r-p) 
\end{eqnarray}

For example, if $f(p,q,r,s)>0$, it means that, against $s_{-i}=aab$, Player $i$'s resulting payoff when playing $s_{i}=a$ is higher than $s_{i}=b$; then $s^*_{i}=a$ for $s_{-i}=aab$; while, if $f(p,q,r,s)<0$, the best response would be $s^*_{i}=b$.

We rearranged the terms in the last equality to be able to simplify the analysis by defining $x\equiv r-p$ and $y\equiv q-s$. By definition, $x>0$ and $y>0$ since our first assumption was $r>p$ and $q>s$. Now, our functions only depend on two parameters $x$ and $y$:

\begin{eqnarray}\label{fgh1}
 & f(x,y)=-2x+y \\
 & g(x,y)=2y-x  \\
 & h(x,y)=y-x 
\end{eqnarray}

The quick analysis gives us four situations with the corresponding best strategy for Player $i$:
\begin{itemize}
    \item[•] $x/2>y$ : play $s^*_i=a$ only when everyone else is playing $b$; play $s^*_i=b$ otherwise.
    \item[•] $x>y>x/2$ : play $s^*_i=a$ when two or more players are playing $b$; play $s^*_i=b$ otherwise.
    \item[•] $2x>y>x$ : play $s^*_i=b$ when two or more players are playing $a$; play $s^*_i=a$ otherwise.
    \item[•] $y>2x$ : play $s^*_i=b$ only when everyone else is playing $a$; play $s^*_i=a$ otherwise.
\end{itemize}
The last two situations are equivalent to the first ones when swapping $a \longleftrightarrow b$. Hence, the study of the first two cases is sufficient to extend the results to the remaining ones.

\mike{Let us make } \textcolor{black}{an additional comment on the analysis performed above. We considered that every player plays the same game with the same payoff parameters $p$, $q$, $r$ and $s$, but \mike{one could have a} situation in which not the same but three different anti-coordination games \mike{are} played\footnote{\textcolor{black}{In the case of  players at the starting and ending points of the ladder, only two different games would be considered.}}. For example, one player plays Game 1 against the player in front, and Game 2 and 3 against the players at both sides. Generically, we would have Game $j$ ($j=1,2,3$) with its associated payoff parameters $p_{j}$, $q_{j}$, $r_{j}$ and $s_{j}$, and the NE condition for the anti-coordination solution $x_{j}\equiv r_{j}-p_{j}>0$ and $y_{j}\equiv q_{j}-s_{j}>0$. This \mike{more general} situation would set up different scenarios depending on the relative values of each $x_{j}$ and $y_{j}$. However, if $x_{1}=x_{2}=x_{3}=x$ and $y_{1}=y_{2}=y_{3}=y$ even when having different $p_{j}$, $q_{j}$, $r_{j}$ and $s_{j}$ for each $j$, the \mike{afore-mentioned} result for the best responses would still apply\footnote{\textcolor{black}{It is easy to see that the corresponding functions to determine the best response would collapse to those in equations ($4$)-($6$).}}.}

\section{Case when $x>y>x/2$}\label{sec:case1}
In this section we will study the case when the payoff parameters obey the relation $x>y>x/2$, which translates into $r-p > q-s >(r-p)/2$.

Let's recall that the rule for the best strategy for each player in this case is: \emph{play $s^*_i=a$ when two or more players are playing $b$, and play $s^*_i=b$ otherwise}\footnote{Another way of remembering the rule for the middle players is: play the opposite of what the majority is playing.}.

\subsection{Ladder}\label{sec:laddercase1}
We want to find the number of different NE solutions for a generic ladder with $2n$ players for this particular case. To do the counting, we will try to simplify the problem: we will use different blocks of 4 players which represent a possible partial solution and the rules to stick each block to the next one. In another words, we will be constructing each solution starting with a block, and then sticking another (allowed) block to it and so on and so forth until the total number of players is reached.  By doing so, the problem of getting the number of NE transforms to a combinatoric one.

The building blocks in this case are:
\[
0 \equiv 
\begin{array}{|c|c|}
    \hline
     \textcolor{red}{a} & \textcolor{blue}{b}  \\ \hline
     \textcolor{blue}{b} & \textcolor{red}{a} \\ \hline
\end{array}
\; \; \; \; \; \; 
1 \equiv 
\begin{array}{|c|c|}
    \hline
     \textcolor{blue}{b} & \textcolor{red}{a}  \\ \hline
     \textcolor{red}{a} & \textcolor{blue}{b} \\ \hline
\end{array}
\; \; \; \; \; \; 
2 \equiv 
\begin{array}{|c|c|}
    \hline
     \textcolor{red}{a} & \textcolor{red}{a}  \\ \hline
     \textcolor{blue}{b} & \textcolor{blue}{b} \\ \hline
\end{array}
\; \; \; \; \; \; 
3 \equiv 
\begin{array}{|c|c|}
    \hline
     \textcolor{blue}{b} & \textcolor{blue}{b}  \\ \hline
     \textcolor{red}{a} & \textcolor{red}{a} \\ \hline
\end{array}
\; \; \; \; \; \; 
4 \equiv
\begin{array}{|c|c|}
    \hline
     \textcolor{blue}{b} & \textcolor{red}{a}  \\ \hline
     \textcolor{blue}{b} & \textcolor{red}{a} \\ \hline
\end{array}
\]

For simplicity, we will label them as blocks $0$, $1$, $2$, $3$ and $4$. As we mentioned, there are some rules to stick each block to the next one, and these also depend whether we want to build the solutions sticking the next block to the right or left of the previous one\footnote{The distinction \mike{in} this particular case doesn't play an important role, but it will in the next subsection.}. We choose sticking them to the left-hand-side of the block.
The rules to get the different solutions are:
\begin{itemize}
    \item No solution can start nor end in blocks $2$ and $3$ because it would imply that the top ending/starting player is choosing $a$ ($b$) against $s_{-i}=ab$ and the bottom ending/starting player would be playing $b$ ($a$) against the same $s_{-i}=ab$, and hence, one of them is not playing a best response.
    \item Before block $0$ there can only be attached blocks $0$ or $1$ or $3$. Hence we can write: $0, 1, 3\rightarrow\textbf{0}$. Similarly for blocks $1$, $2$ and $3$, we can argue: $0, 1, 2\rightarrow\textbf{1}$ ; $0, 3 \rightarrow\textbf{2}$ and $1, 2\rightarrow\textbf{3}$.
    \item The block $4$ can only be attached to itself: $4\rightarrow \textbf{4}$ because attaching it to a different one would imply that one player would be playing $s_i=a$ ($s_i=b$) against $s_{-i}=aab$ ($s_{-i}=abb$) and that would not be a best response. In addition, the last two players at the edge of the ladder cannot be playing $|a|a|$ \footnote{For simplicity, we wrote those blocks horizontally, but they should be read vertically: $|a|b|$ means that on the top row the player chooses to play $a$, while the one in front plays $b$.} because they would be playing $s_i=a$ against $s_{-i}=ab$, which, again, would not be a best response.
\end{itemize}

\textcolor{black}{All these rules can be more easily understood by formulating them as rules in the graph. Each player/node that chooses to play $a$ is surrounded by red lines, while if $b$, the lines are blue. Each time a red line meets a blue one, the full line turns to green. The rule in the graph to connect each block to get all the NE is: each node must have \emph{at least two green} lines\footnote{\textcolor{black}{There is an exception to this when using block $4$ with odd $n$, in which the starting and ending nodes have only one green line.}}. An example of that will be explained and shown below.}

Now, obtaining all the possible combinations given those rules becomes an easy task, but first let's focus on the relation between the total number of players $2n$ and its relation to the total number of 4-player blocks needed.

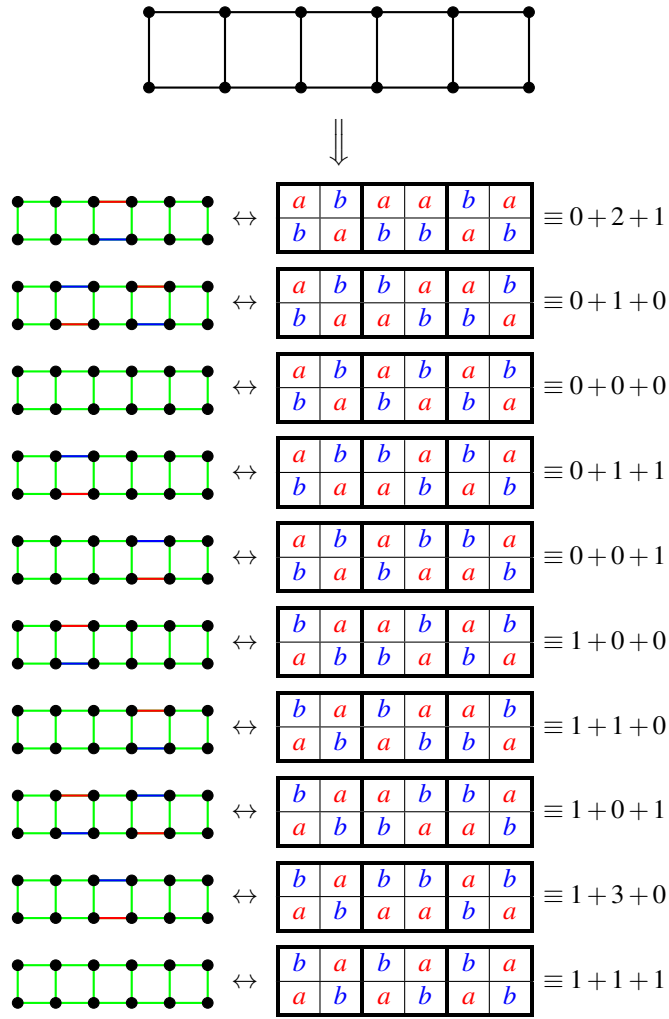
\begin{figure}
\centering
\begin{minipage}[c]{1\linewidth}
\[
\begin{tikzpicture}
\draw[black, thick] (0,-0.5) -- (5,-0.5);
\draw[black, thick] (0,0.5) -- (5,0.5);
\draw[black, thick] (0,0.5) -- (0,-0.5) ;
\draw[black, thick] (1,0.5) -- (1,-0.5) ;
\draw[black, thick] (2,0.5) -- (2,-0.5) ;
\draw[black, thick] (3,0.5) -- (3,-0.5);
\draw[black, thick] (4,0.5) -- (4,-0.5);
\draw[black, thick] (5,0.5) -- (5,-0.5);
\filldraw[color=black] (0,-0.5) circle (2pt);
\filldraw[color=black] (1,-0.5) circle (2pt);
\filldraw[color=black] (2,-0.5) circle (2pt);
\filldraw[color=black] (3,-0.5) circle (2pt);
\filldraw[color=black] (4,-0.5) circle (2pt);
\filldraw[color=black] (5,-0.5) circle (2pt);
\filldraw[color=black] (0,0.5) circle (2pt);
\filldraw[color=black] (1,0.5) circle (2pt);
\filldraw[color=black] (2,0.5) circle (2pt);
\filldraw[color=black] (3,0.5) circle (2pt);
\filldraw[color=black] (4,0.5) circle (2pt);
\filldraw[color=black] (5,0.5) circle (2pt);
\end{tikzpicture}
\vspace{3pt}
  \]
 \[ \Big\Downarrow \]
\end{minipage} \vspace{3pt}

\begin{minipage}{1\linewidth}
\[
\begin{tikzpicture}[anchor=base, baseline=-1.5pt]
\draw[green, thick] (0,-0.25) -- (2.5,-0.25);
\draw[green, thick] (0,0.25) -- (2.5,0.25);
\draw[green, thick] (0,0.25) -- (0,-0.25) ;
\draw[green, thick] (0.5,0.25) -- (0.5,-0.25) ;
\draw[green, thick] (1,0.25) -- (1,-0.25) ;
\draw[green, thick] (1.5,0.25) -- (1.5,-0.25);
\draw[green, thick] (2,0.25) -- (2,-0.25);
\draw[green, thick] (2.5,0.25) -- (2.5,-0.25);
\draw[red, thick] (1,0.25) -- (1.5,0.25);
\draw[blue, thick] (1,-0.25) -- (1.5,-0.25);
\filldraw[color=black] (0,-0.25) circle (2pt);
\filldraw[color=black] (0.5,-0.25) circle (2pt);
\filldraw[color=black] (1,-0.25) circle (2pt);
\filldraw[color=black] (1.5,-0.25) circle (2pt);
\filldraw[color=black] (2,-0.25) circle (2pt);
\filldraw[color=black] (2.5,-0.25) circle (2pt);
\filldraw[color=black] (0,0.25) circle (2pt);
\filldraw[color=black] (0.5,0.25) circle (2pt);
\filldraw[color=black] (1,0.25) circle (2pt);
\filldraw[color=black] (1.5,0.25) circle (2pt);
\filldraw[color=black] (2,0.25) circle (2pt);
\filldraw[color=black] (2.5,0.25) circle (2pt);
\end{tikzpicture}
\; \;  \leftrightarrow \; \; 
\begin{tabu}{|[1.3pt]c|c|[1.3pt]c|c|[1.3pt]c|c|[1.3pt]}
    \tabucline[1.3pt]{-} 
     \textcolor{red}{a} & \textcolor{blue}{b} & \textcolor{red}{a} & \textcolor{red}{a} & \textcolor{blue}{b} & \textcolor{red}{a} \\ \hline 
     \textcolor{blue}{b} & \textcolor{red}{a} & \textcolor{blue}{b} & \textcolor{blue}{b} & \textcolor{red}{a} & \textcolor{blue}{b}\\ \tabucline[1.3pt]{-}
\end{tabu} 
\:  \equiv 0+2+1 
\]
\end{minipage}\vspace{3pt}

\begin{minipage}{1\linewidth}
\[
\begin{tikzpicture}[anchor=base, baseline=-1.5pt]
\draw[green, thick] (0,-0.25) -- (2.5,-0.25);
\draw[green, thick] (0,0.25) -- (2.5,0.25);
\draw[green, thick] (0,0.25) -- (0,-0.25) ;
\draw[green, thick] (0.5,0.25) -- (0.5,-0.25) ;
\draw[green, thick] (1,0.25) -- (1,-0.25) ;
\draw[green, thick] (1.5,0.25) -- (1.5,-0.25);
\draw[green, thick] (2,0.25) -- (2,-0.25);
\draw[green, thick] (2.5,0.25) -- (2.5,-0.25);
\draw[blue, thick] (0.5,0.25) -- (1,0.25);
\draw[red, thick] (0.5,-0.25) -- (1,-0.25);
\draw[red, thick] (1.5,0.25) -- (2,0.25);
\draw[blue, thick] (1.5,-0.25) -- (2,-0.25);
\filldraw[color=black] (0,-0.25) circle (2pt);
\filldraw[color=black] (0.5,-0.25) circle (2pt);
\filldraw[color=black] (1,-0.25) circle (2pt);
\filldraw[color=black] (1.5,-0.25) circle (2pt);
\filldraw[color=black] (2,-0.25) circle (2pt);
\filldraw[color=black] (2.5,-0.25) circle (2pt);
\filldraw[color=black] (0,0.25) circle (2pt);
\filldraw[color=black] (0.5,0.25) circle (2pt);
\filldraw[color=black] (1,0.25) circle (2pt);
\filldraw[color=black] (1.5,0.25) circle (2pt);
\filldraw[color=black] (2,0.25) circle (2pt);
\filldraw[color=black] (2.5,0.25) circle (2pt);
\end{tikzpicture}
\; \;  \leftrightarrow \; \; 
\begin{tabu}{|[1.3pt]c|c|[1.3pt]c|c|[1.3pt]c|c|[1.3pt]}
    \tabucline[1.3pt]{-} 
     \textcolor{red}{a} & \textcolor{blue}{b} & \textcolor{blue}{b} & \textcolor{red}{a} & \textcolor{red}{a} & \textcolor{blue}{b} \\ \hline 
     \textcolor{blue}{b} & \textcolor{red}{a} & \textcolor{red}{a} & \textcolor{blue}{b} & \textcolor{blue}{b} & \textcolor{red}{a}\\ \tabucline[1.3pt]{-}
\end{tabu} 
\:  \equiv 0+1+0
\]
\end{minipage}\vspace{3pt}

\begin{minipage}{1\linewidth}
\[
\begin{tikzpicture}[anchor=base, baseline=-1.5pt]
\draw[green, thick] (0,-0.25) -- (2.5,-0.25);
\draw[green, thick] (0,0.25) -- (2.5,0.25);
\draw[green, thick] (0,0.25) -- (0,-0.25) ;
\draw[green, thick] (0.5,0.25) -- (0.5,-0.25) ;
\draw[green, thick] (1,0.25) -- (1,-0.25) ;
\draw[green, thick] (1.5,0.25) -- (1.5,-0.25);
\draw[green, thick] (2,0.25) -- (2,-0.25);
\draw[green, thick] (2.5,0.25) -- (2.5,-0.25);
\filldraw[color=black] (0,-0.25) circle (2pt);
\filldraw[color=black] (0.5,-0.25) circle (2pt);
\filldraw[color=black] (1,-0.25) circle (2pt);
\filldraw[color=black] (1.5,-0.25) circle (2pt);
\filldraw[color=black] (2,-0.25) circle (2pt);
\filldraw[color=black] (2.5,-0.25) circle (2pt);
\filldraw[color=black] (0,0.25) circle (2pt);
\filldraw[color=black] (0.5,0.25) circle (2pt);
\filldraw[color=black] (1,0.25) circle (2pt);
\filldraw[color=black] (1.5,0.25) circle (2pt);
\filldraw[color=black] (2,0.25) circle (2pt);
\filldraw[color=black] (2.5,0.25) circle (2pt);
\end{tikzpicture}
\; \;  \leftrightarrow \; \; 
\begin{tabu}{|[1.3pt]c|c|[1.3pt]c|c|[1.3pt]c|c|[1.3pt]}
    \tabucline[1.3pt]{-} 
     \textcolor{red}{a} & \textcolor{blue}{b} & \textcolor{red}{a} & \textcolor{blue}{b} & \textcolor{red}{a} & \textcolor{blue}{b} \\ \hline 
     \textcolor{blue}{b} & \textcolor{red}{a} & \textcolor{blue}{b} & \textcolor{red}{a} & \textcolor{blue}{b} & \textcolor{red}{a}\\ \tabucline[1.3pt]{-}
\end{tabu} 
\:  \equiv 0+0+0
\]
\end{minipage}\vspace{3pt}

\begin{minipage}{1\linewidth}
\[
\begin{tikzpicture}[anchor=base, baseline=-1.5pt]
\draw[green, thick] (0,-0.25) -- (2.5,-0.25);
\draw[green, thick] (0,0.25) -- (2.5,0.25);
\draw[green, thick] (0,0.25) -- (0,-0.25) ;
\draw[green, thick] (0.5,0.25) -- (0.5,-0.25) ;
\draw[green, thick] (1,0.25) -- (1,-0.25) ;
\draw[green, thick] (1.5,0.25) -- (1.5,-0.25);
\draw[green, thick] (2,0.25) -- (2,-0.25);
\draw[green, thick] (2.5,0.25) -- (2.5,-0.25);
\draw[blue, thick] (0.5,0.25) -- (1,0.25);
\draw[red, thick] (0.5,-0.25) -- (1,-0.25);
\filldraw[color=black] (0,-0.25) circle (2pt);
\filldraw[color=black] (0.5,-0.25) circle (2pt);
\filldraw[color=black] (1,-0.25) circle (2pt);
\filldraw[color=black] (1.5,-0.25) circle (2pt);
\filldraw[color=black] (2,-0.25) circle (2pt);
\filldraw[color=black] (2.5,-0.25) circle (2pt);
\filldraw[color=black] (0,0.25) circle (2pt);
\filldraw[color=black] (0.5,0.25) circle (2pt);
\filldraw[color=black] (1,0.25) circle (2pt);
\filldraw[color=black] (1.5,0.25) circle (2pt);
\filldraw[color=black] (2,0.25) circle (2pt);
\filldraw[color=black] (2.5,0.25) circle (2pt);
\end{tikzpicture}
\; \;  \leftrightarrow \; \; 
\begin{tabu}{|[1.3pt]c|c|[1.3pt]c|c|[1.3pt]c|c|[1.3pt]}
    \tabucline[1.3pt]{-} 
     \textcolor{red}{a} & \textcolor{blue}{b} & \textcolor{blue}{b} & \textcolor{red}{a} & \textcolor{blue}{b} & \textcolor{red}{a} \\ \hline 
     \textcolor{blue}{b} & \textcolor{red}{a} & \textcolor{red}{a} & \textcolor{blue}{b} & \textcolor{red}{a} & \textcolor{blue}{b}\\ \tabucline[1.3pt]{-}
\end{tabu} 
\:  \equiv 0+1+1
\]
\end{minipage}\vspace{3pt}

\begin{minipage}{1\linewidth}
\[
\begin{tikzpicture}[anchor=base, baseline=-1.5pt]
\draw[green, thick] (0,-0.25) -- (2.5,-0.25);
\draw[green, thick] (0,0.25) -- (2.5,0.25);
\draw[green, thick] (0,0.25) -- (0,-0.25) ;
\draw[green, thick] (0.5,0.25) -- (0.5,-0.25) ;
\draw[green, thick] (1,0.25) -- (1,-0.25) ;
\draw[green, thick] (1.5,0.25) -- (1.5,-0.25);
\draw[green, thick] (2,0.25) -- (2,-0.25);
\draw[green, thick] (2.5,0.25) -- (2.5,-0.25);
\draw[blue, thick] (1.5,0.25) -- (2,0.25);
\draw[red, thick] (1.5,-0.25) -- (2,-0.25);
\filldraw[color=black] (0,-0.25) circle (2pt);
\filldraw[color=black] (0.5,-0.25) circle (2pt);
\filldraw[color=black] (1,-0.25) circle (2pt);
\filldraw[color=black] (1.5,-0.25) circle (2pt);
\filldraw[color=black] (2,-0.25) circle (2pt);
\filldraw[color=black] (2.5,-0.25) circle (2pt);
\filldraw[color=black] (0,0.25) circle (2pt);
\filldraw[color=black] (0.5,0.25) circle (2pt);
\filldraw[color=black] (1,0.25) circle (2pt);
\filldraw[color=black] (1.5,0.25) circle (2pt);
\filldraw[color=black] (2,0.25) circle (2pt);
\filldraw[color=black] (2.5,0.25) circle (2pt);
\end{tikzpicture}
\; \;  \leftrightarrow \; \; 
\begin{tabu}{|[1.3pt]c|c|[1.3pt]c|c|[1.3pt]c|c|[1.3pt]}
    \tabucline[1.3pt]{-} 
     \textcolor{red}{a} & \textcolor{blue}{b} & \textcolor{red}{a} & \textcolor{blue}{b} & \textcolor{blue}{b} & \textcolor{red}{a} \\ \hline 
     \textcolor{blue}{b} & \textcolor{red}{a} & \textcolor{blue}{b} & \textcolor{red}{a} & \textcolor{red}{a} & \textcolor{blue}{b}\\ \tabucline[1.3pt]{-}
\end{tabu} 
\:  \equiv 0+0+1
\]
\end{minipage}\vspace{3pt}

\begin{minipage}{1\linewidth}
\[
\begin{tikzpicture}[anchor=base, baseline=-1.5pt]
\draw[green, thick] (0,-0.25) -- (2.5,-0.25);
\draw[green, thick] (0,0.25) -- (2.5,0.25);
\draw[green, thick] (0,0.25) -- (0,-0.25) ;
\draw[green, thick] (0.5,0.25) -- (0.5,-0.25) ;
\draw[green, thick] (1,0.25) -- (1,-0.25) ;
\draw[green, thick] (1.5,0.25) -- (1.5,-0.25);
\draw[green, thick] (2,0.25) -- (2,-0.25);
\draw[green, thick] (2.5,0.25) -- (2.5,-0.25);
\draw[red, thick] (0.5,0.25) -- (1,0.25);
\draw[blue, thick] (0.5,-0.25) -- (1,-0.25);
\filldraw[color=black] (0,-0.25) circle (2pt);
\filldraw[color=black] (0.5,-0.25) circle (2pt);
\filldraw[color=black] (1,-0.25) circle (2pt);
\filldraw[color=black] (1.5,-0.25) circle (2pt);
\filldraw[color=black] (2,-0.25) circle (2pt);
\filldraw[color=black] (2.5,-0.25) circle (2pt);
\filldraw[color=black] (0,0.25) circle (2pt);
\filldraw[color=black] (0.5,0.25) circle (2pt);
\filldraw[color=black] (1,0.25) circle (2pt);
\filldraw[color=black] (1.5,0.25) circle (2pt);
\filldraw[color=black] (2,0.25) circle (2pt);
\filldraw[color=black] (2.5,0.25) circle (2pt);
\end{tikzpicture}
\; \;  \leftrightarrow \; \; 
\begin{tabu}{|[1.3pt]c|c|[1.3pt]c|c|[1.3pt]c|c|[1.3pt]}
    \tabucline[1.3pt]{-} 
     \textcolor{blue}{b} & \textcolor{red}{a} & \textcolor{red}{a} & \textcolor{blue}{b} & \textcolor{red}{a} & \textcolor{blue}{b} \\ \hline 
     \textcolor{red}{a} & \textcolor{blue}{b} & \textcolor{blue}{b} & \textcolor{red}{a} & \textcolor{blue}{b} & \textcolor{red}{a}\\ \tabucline[1.3pt]{-}
\end{tabu} 
\:  \equiv 1+0+0
\]
\end{minipage}\vspace{3pt}

\begin{minipage}{1\linewidth}
\[
\begin{tikzpicture}[anchor=base, baseline=-1.5pt]
\draw[green, thick] (0,-0.25) -- (2.5,-0.25);
\draw[green, thick] (0,0.25) -- (2.5,0.25);
\draw[green, thick] (0,0.25) -- (0,-0.25) ;
\draw[green, thick] (0.5,0.25) -- (0.5,-0.25) ;
\draw[green, thick] (1,0.25) -- (1,-0.25) ;
\draw[green, thick] (1.5,0.25) -- (1.5,-0.25);
\draw[green, thick] (2,0.25) -- (2,-0.25);
\draw[green, thick] (2.5,0.25) -- (2.5,-0.25);
\draw[red, thick] (1.5,0.25) -- (2,0.25);
\draw[blue, thick] (1.5,-0.25) -- (2,-0.25);
\filldraw[color=black] (0,-0.25) circle (2pt);
\filldraw[color=black] (0.5,-0.25) circle (2pt);
\filldraw[color=black] (1,-0.25) circle (2pt);
\filldraw[color=black] (1.5,-0.25) circle (2pt);
\filldraw[color=black] (2,-0.25) circle (2pt);
\filldraw[color=black] (2.5,-0.25) circle (2pt);
\filldraw[color=black] (0,0.25) circle (2pt);
\filldraw[color=black] (0.5,0.25) circle (2pt);
\filldraw[color=black] (1,0.25) circle (2pt);
\filldraw[color=black] (1.5,0.25) circle (2pt);
\filldraw[color=black] (2,0.25) circle (2pt);
\filldraw[color=black] (2.5,0.25) circle (2pt);
\end{tikzpicture}
\; \;  \leftrightarrow \; \; 
\begin{tabu}{|[1.3pt]c|c|[1.3pt]c|c|[1.3pt]c|c|[1.3pt]}
    \tabucline[1.3pt]{-} 
     \textcolor{blue}{b} & \textcolor{red}{a} & \textcolor{blue}{b} & \textcolor{red}{a} & \textcolor{red}{a} & \textcolor{blue}{b} \\ \hline 
     \textcolor{red}{a} & \textcolor{blue}{b} & \textcolor{red}{a} & \textcolor{blue}{b} & \textcolor{blue}{b} & \textcolor{red}{a}\\ \tabucline[1.3pt]{-}
\end{tabu} 
\:  \equiv 1+1+0
\]
\end{minipage}\vspace{3pt}

\begin{minipage}{1\linewidth}
\[
\begin{tikzpicture}[anchor=base, baseline=-1.5pt]
\draw[green, thick] (0,-0.25) -- (2.5,-0.25);
\draw[green, thick] (0,0.25) -- (2.5,0.25);
\draw[green, thick] (0,0.25) -- (0,-0.25) ;
\draw[green, thick] (0.5,0.25) -- (0.5,-0.25) ;
\draw[green, thick] (1,0.25) -- (1,-0.25) ;
\draw[green, thick] (1.5,0.25) -- (1.5,-0.25);
\draw[green, thick] (2,0.25) -- (2,-0.25);
\draw[green, thick] (2.5,0.25) -- (2.5,-0.25);
\draw[red, thick] (0.5,0.25) -- (1,0.25);
\draw[red, thick] (1.5,-0.25) -- (2,-0.25);
\draw[blue, thick] (0.5,-0.25) -- (1,-0.25);
\draw[blue, thick] (1.5,0.25) -- (2,0.25);
\filldraw[color=black] (0,-0.25) circle (2pt);
\filldraw[color=black] (0.5,-0.25) circle (2pt);
\filldraw[color=black] (1,-0.25) circle (2pt);
\filldraw[color=black] (1.5,-0.25) circle (2pt);
\filldraw[color=black] (2,-0.25) circle (2pt);
\filldraw[color=black] (2.5,-0.25) circle (2pt);
\filldraw[color=black] (0,0.25) circle (2pt);
\filldraw[color=black] (0.5,0.25) circle (2pt);
\filldraw[color=black] (1,0.25) circle (2pt);
\filldraw[color=black] (1.5,0.25) circle (2pt);
\filldraw[color=black] (2,0.25) circle (2pt);
\filldraw[color=black] (2.5,0.25) circle (2pt);
\end{tikzpicture}
\; \;  \leftrightarrow \; \; 
\begin{tabu}{|[1.3pt]c|c|[1.3pt]c|c|[1.3pt]c|c|[1.3pt]}
    \tabucline[1.3pt]{-} 
     \textcolor{blue}{b} & \textcolor{red}{a} & \textcolor{red}{a} & \textcolor{blue}{b} & \textcolor{blue}{b} & \textcolor{red}{a} \\ \hline 
     \textcolor{red}{a} & \textcolor{blue}{b} & \textcolor{blue}{b} & \textcolor{red}{a} & \textcolor{red}{a} & \textcolor{blue}{b}\\ \tabucline[1.3pt]{-}
\end{tabu} 
\:  \equiv 1+0+1
\]
\end{minipage}\vspace{3pt}

\begin{minipage}{1\linewidth}
\[
\begin{tikzpicture}[anchor=base, baseline=-1.5pt]
\draw[green, thick] (0,-0.25) -- (2.5,-0.25);
\draw[green, thick] (0,0.25) -- (2.5,0.25);
\draw[green, thick] (0,0.25) -- (0,-0.25) ;
\draw[green, thick] (0.5,0.25) -- (0.5,-0.25) ;
\draw[green, thick] (1,0.25) -- (1,-0.25) ;
\draw[green, thick] (1.5,0.25) -- (1.5,-0.25);
\draw[green, thick] (2,0.25) -- (2,-0.25);
\draw[green, thick] (2.5,0.25) -- (2.5,-0.25);
\draw[red, thick] (1,-0.25) -- (1.5,-0.25);
\draw[blue, thick] (1,0.25) -- (1.5,0.25);
\filldraw[color=black] (0,-0.25) circle (2pt);
\filldraw[color=black] (0.5,-0.25) circle (2pt);
\filldraw[color=black] (1,-0.25) circle (2pt);
\filldraw[color=black] (1.5,-0.25) circle (2pt);
\filldraw[color=black] (2,-0.25) circle (2pt);
\filldraw[color=black] (2.5,-0.25) circle (2pt);
\filldraw[color=black] (0,0.25) circle (2pt);
\filldraw[color=black] (0.5,0.25) circle (2pt);
\filldraw[color=black] (1,0.25) circle (2pt);
\filldraw[color=black] (1.5,0.25) circle (2pt);
\filldraw[color=black] (2,0.25) circle (2pt);
\filldraw[color=black] (2.5,0.25) circle (2pt);
\end{tikzpicture}
\; \;  \leftrightarrow \; \; 
\begin{tabu}{|[1.3pt]c|c|[1.3pt]c|c|[1.3pt]c|c|[1.3pt]}
    \tabucline[1.3pt]{-} 
     \textcolor{blue}{b} & \textcolor{red}{a} & \textcolor{blue}{b} & \textcolor{blue}{b} & \textcolor{red}{a} & \textcolor{blue}{b} \\ \hline 
     \textcolor{red}{a} & \textcolor{blue}{b} & \textcolor{red}{a} & \textcolor{red}{a} & \textcolor{blue}{b} & \textcolor{red}{a}\\ \tabucline[1.3pt]{-}
\end{tabu} 
\:  \equiv 1+3+0
\]
\end{minipage}\vspace{3pt}

\begin{minipage}{1\linewidth}
\[
\begin{tikzpicture}[anchor=base, baseline=-1.5pt]
\draw[green, thick] (0,-0.25) -- (2.5,-0.25);
\draw[green, thick] (0,0.25) -- (2.5,0.25);
\draw[green, thick] (0,0.25) -- (0,-0.25) ;
\draw[green, thick] (0.5,0.25) -- (0.5,-0.25) ;
\draw[green, thick] (1,0.25) -- (1,-0.25) ;
\draw[green, thick] (1.5,0.25) -- (1.5,-0.25);
\draw[green, thick] (2,0.25) -- (2,-0.25);
\draw[green, thick] (2.5,0.25) -- (2.5,-0.25);
\filldraw[color=black] (0,-0.25) circle (2pt);
\filldraw[color=black] (0.5,-0.25) circle (2pt);
\filldraw[color=black] (1,-0.25) circle (2pt);
\filldraw[color=black] (1.5,-0.25) circle (2pt);
\filldraw[color=black] (2,-0.25) circle (2pt);
\filldraw[color=black] (2.5,-0.25) circle (2pt);
\filldraw[color=black] (0,0.25) circle (2pt);
\filldraw[color=black] (0.5,0.25) circle (2pt);
\filldraw[color=black] (1,0.25) circle (2pt);
\filldraw[color=black] (1.5,0.25) circle (2pt);
\filldraw[color=black] (2,0.25) circle (2pt);
\filldraw[color=black] (2.5,0.25) circle (2pt);
\end{tikzpicture}
\; \;  \leftrightarrow \; \; 
\begin{tabu}{|[1.3pt]c|c|[1.3pt]c|c|[1.3pt]c|c|[1.3pt]}
    \tabucline[1.3pt]{-} 
     \textcolor{blue}{b} & \textcolor{red}{a} & \textcolor{blue}{b} & \textcolor{red}{a} & \textcolor{blue}{b} & \textcolor{red}{a} \\ \hline 
     \textcolor{red}{a} & \textcolor{blue}{b} & \textcolor{red}{a} & \textcolor{blue}{b} & \textcolor{red}{a} & \textcolor{blue}{b}\\ \tabucline[1.3pt]{-}
\end{tabu} 
\:  \equiv 1+1+1
\]
\end{minipage}
\caption{All 10 possible solutions for the ladder in the case when $x>y>x/2$ with $2n=12$ players decomposed in the elementary blocks \textcolor{black}{and the graphic representation}. Each player plays a game with their neighbour (to the left, right and on top/bottom). \textcolor{black}{The green line denotes a matching of a red and blue line, that is, of two players, one of \mike{whom chooses} $a$ and the other $b$. The rule states that every node must be surrounded by at least two green lines.}}
\label{fig:example12}
\end{figure}

In the case when $2n$ is a multiple of 4, i.e. $2n=4k$ (equivalently, $n$ is even), in which $k$ stands for the number of blocks of 4 players, we would need to find all the possible combinations of $k$ blocks using $0$, $1$, $2$ and $3$ \footnote{We ignore the solution with block $4$ because the unique combination of $4444..44$ is not allowed due to the last two players choosing $|a|a|$.} given the rules above. In figure \ref{fig:example12} we show all the possible solutions for 12 players ($k=4$ blocks) and the decomposition in the elemental blocks.

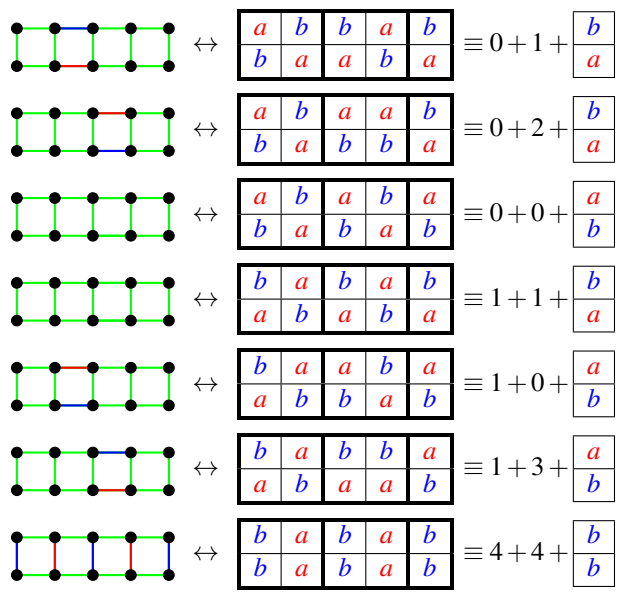
\begin{figure}
\centering
\begin{minipage}{1\linewidth}
\[
\begin{tikzpicture}[anchor=base, baseline=-1.5pt]
\draw[green, thick] (0,-0.25) -- (2,-0.25);
\draw[green, thick] (0,0.25) -- (2,0.25);
\draw[green, thick] (0,0.25) -- (0,-0.25) ;
\draw[green, thick] (0.5,0.25) -- (0.5,-0.25) ;
\draw[green, thick] (1,0.25) -- (1,-0.25) ;
\draw[green, thick] (1.5,0.25) -- (1.5,-0.25);
\draw[green, thick] (2,0.25) -- (2,-0.25);
\draw[blue, thick] (0.5,0.25) -- (1,0.25);
\draw[red, thick] (0.5,-0.25) -- (1,-0.25);
\filldraw[color=black] (0,-0.25) circle (2pt);
\filldraw[color=black] (0.5,-0.25) circle (2pt);
\filldraw[color=black] (1,-0.25) circle (2pt);
\filldraw[color=black] (1.5,-0.25) circle (2pt);
\filldraw[color=black] (2,-0.25) circle (2pt);
\filldraw[color=black] (0,0.25) circle (2pt);
\filldraw[color=black] (0.5,0.25) circle (2pt);
\filldraw[color=black] (1,0.25) circle (2pt);
\filldraw[color=black] (1.5,0.25) circle (2pt);
\filldraw[color=black] (2,0.25) circle (2pt);
\end{tikzpicture}
\; \; \leftrightarrow \; \;
\begin{tabu}{|[1.3pt]c|c|[1.3pt]c|c|[1.3pt]c|[1.3pt]}
    \tabucline[1.3pt]{-} 
     \textcolor{red}{a} & \textcolor{blue}{b} & \textcolor{blue}{b} & \textcolor{red}{a} & \textcolor{blue}{b} \\ \hline 
     \textcolor{blue}{b} & \textcolor{red}{a} & \textcolor{red}{a} & \textcolor{blue}{b} & \textcolor{red}{a}\\ \tabucline[1.3pt]{-}
\end{tabu} 
\:  \equiv 0+1+
\begin{array}{|c|}
    \hline
     \textcolor{blue}{b} \\ \hline
     \textcolor{red}{a} \\ \hline
\end{array}
\]
\end{minipage}\vspace{3pt}

\begin{minipage}{1\linewidth}
\[
\begin{tikzpicture}[anchor=base, baseline=-1.5pt]
\draw[green, thick] (0,-0.25) -- (2,-0.25);
\draw[green, thick] (0,0.25) -- (2,0.25);
\draw[green, thick] (0,0.25) -- (0,-0.25) ;
\draw[green, thick] (0.5,0.25) -- (0.5,-0.25) ;
\draw[green, thick] (1,0.25) -- (1,-0.25) ;
\draw[green, thick] (1.5,0.25) -- (1.5,-0.25);
\draw[green, thick] (2,0.25) -- (2,-0.25);
\draw[red, thick] (1,0.25) -- (1.5,0.25);
\draw[blue, thick] (1,-0.25) -- (1.5,-0.25);
\filldraw[color=black] (0,-0.25) circle (2pt);
\filldraw[color=black] (0.5,-0.25) circle (2pt);
\filldraw[color=black] (1,-0.25) circle (2pt);
\filldraw[color=black] (1.5,-0.25) circle (2pt);
\filldraw[color=black] (2,-0.25) circle (2pt);
\filldraw[color=black] (0,0.25) circle (2pt);
\filldraw[color=black] (0.5,0.25) circle (2pt);
\filldraw[color=black] (1,0.25) circle (2pt);
\filldraw[color=black] (1.5,0.25) circle (2pt);
\filldraw[color=black] (2,0.25) circle (2pt);
\end{tikzpicture}
\; \; \leftrightarrow \; \;
\begin{tabu}{|[1.3pt]c|c|[1.3pt]c|c|[1.3pt]c|[1.3pt]}
    \tabucline[1.3pt]{-}
     \textcolor{red}{a} & \textcolor{blue}{b} & \textcolor{red}{a} & \textcolor{red}{a} & \textcolor{blue}{b} \\ \hline 
     \textcolor{blue}{b} & \textcolor{red}{a} & \textcolor{blue}{b} & \textcolor{blue}{b} & \textcolor{red}{a}\\ \tabucline[1.3pt]{-}
\end{tabu} 
\:  \equiv 0+2+
\begin{array}{|c|}
    \hline
     \textcolor{blue}{b} \\ \hline
     \textcolor{red}{a} \\ \hline
\end{array}
\]
\end{minipage}\vspace{3pt}

\begin{minipage}{1\linewidth}
\[
\begin{tikzpicture}[anchor=base, baseline=-1.5pt]
\draw[green, thick] (0,-0.25) -- (2,-0.25);
\draw[green, thick] (0,0.25) -- (2,0.25);
\draw[green, thick] (0,0.25) -- (0,-0.25) ;
\draw[green, thick] (0.5,0.25) -- (0.5,-0.25) ;
\draw[green, thick] (1,0.25) -- (1,-0.25) ;
\draw[green, thick] (1.5,0.25) -- (1.5,-0.25);
\draw[green, thick] (2,0.25) -- (2,-0.25);
\draw[green, thick] (1,0.25) -- (1.5,0.25);
\draw[green, thick] (1,-0.25) -- (1.5,-0.25);
\filldraw[color=black] (0,-0.25) circle (2pt);
\filldraw[color=black] (0.5,-0.25) circle (2pt);
\filldraw[color=black] (1,-0.25) circle (2pt);
\filldraw[color=black] (1.5,-0.25) circle (2pt);
\filldraw[color=black] (2,-0.25) circle (2pt);
\filldraw[color=black] (0,0.25) circle (2pt);
\filldraw[color=black] (0.5,0.25) circle (2pt);
\filldraw[color=black] (1,0.25) circle (2pt);
\filldraw[color=black] (1.5,0.25) circle (2pt);
\filldraw[color=black] (2,0.25) circle (2pt);
\end{tikzpicture}
\; \; \leftrightarrow \; \;
\begin{tabu}{|[1.3pt]c|c|[1.3pt]c|c|[1.3pt]c|[1.3pt]}
    \tabucline[1.3pt]{-}
     \textcolor{red}{a} & \textcolor{blue}{b} & \textcolor{red}{a} & \textcolor{blue}{b} & \textcolor{red}{a} \\ \hline 
     \textcolor{blue}{b} & \textcolor{red}{a} & \textcolor{blue}{b} & \textcolor{red}{a} & \textcolor{blue}{b}\\ \tabucline[1.3pt]{-}
\end{tabu} 
\: \equiv 0+0+
\begin{array}{|c|}
    \hline
     \textcolor{red}{a} \\ \hline
     \textcolor{blue}{b} \\ \hline
\end{array}
\]
\end{minipage}\vspace{3pt}

\begin{minipage}{1\linewidth}
\[
\begin{tikzpicture}[anchor=base, baseline=-1.5pt]
\draw[green, thick] (0,-0.25) -- (2,-0.25);
\draw[green, thick] (0,0.25) -- (2,0.25);
\draw[green, thick] (0,0.25) -- (0,-0.25) ;
\draw[green, thick] (0.5,0.25) -- (0.5,-0.25) ;
\draw[green, thick] (1,0.25) -- (1,-0.25) ;
\draw[green, thick] (1.5,0.25) -- (1.5,-0.25);
\draw[green, thick] (2,0.25) -- (2,-0.25);
\draw[green, thick] (1,0.25) -- (1.5,0.25);
\draw[green, thick] (1,-0.25) -- (1.5,-0.25);
\filldraw[color=black] (0,-0.25) circle (2pt);
\filldraw[color=black] (0.5,-0.25) circle (2pt);
\filldraw[color=black] (1,-0.25) circle (2pt);
\filldraw[color=black] (1.5,-0.25) circle (2pt);
\filldraw[color=black] (2,-0.25) circle (2pt);
\filldraw[color=black] (0,0.25) circle (2pt);
\filldraw[color=black] (0.5,0.25) circle (2pt);
\filldraw[color=black] (1,0.25) circle (2pt);
\filldraw[color=black] (1.5,0.25) circle (2pt);
\filldraw[color=black] (2,0.25) circle (2pt);
\end{tikzpicture}
\; \; \leftrightarrow \; \;
\begin{tabu}{|[1.3pt]c|c|[1.3pt]c|c|[1.3pt]c|[1.3pt]}
    \tabucline[1.3pt]{-}
     \textcolor{blue}{b} & \textcolor{red}{a} & \textcolor{blue}{b} & \textcolor{red}{a} & \textcolor{blue}{b} \\ \hline 
     \textcolor{red}{a} & \textcolor{blue}{b} & \textcolor{red}{a} & \textcolor{blue}{b} & \textcolor{red}{a}\\ \tabucline[1.3pt]{-}
\end{tabu} 
\:  \equiv 1+1+
\begin{array}{|c|}
    \hline
     \textcolor{blue}{b} \\ \hline
     \textcolor{red}{a} \\ \hline
\end{array}
\]
\end{minipage}\vspace{3pt}

\begin{minipage}{1\linewidth}
\[
\begin{tikzpicture}[anchor=base, baseline=-1.5pt]
\draw[green, thick] (0,-0.25) -- (2,-0.25);
\draw[green, thick] (0,0.25) -- (2,0.25);
\draw[green, thick] (0,0.25) -- (0,-0.25) ;
\draw[green, thick] (0.5,0.25) -- (0.5,-0.25) ;
\draw[green, thick] (1,0.25) -- (1,-0.25) ;
\draw[green, thick] (1.5,0.25) -- (1.5,-0.25);
\draw[green, thick] (2,0.25) -- (2,-0.25);
\draw[red, thick] (0.5, 0.25) -- (1, 0.25);
\draw[blue, thick] (0.5,-0.25) -- (1,-0.25);
\filldraw[color=black] (0,-0.25) circle (2pt);
\filldraw[color=black] (0.5,-0.25) circle (2pt);
\filldraw[color=black] (1,-0.25) circle (2pt);
\filldraw[color=black] (1.5,-0.25) circle (2pt);
\filldraw[color=black] (2,-0.25) circle (2pt);
\filldraw[color=black] (0,0.25) circle (2pt);
\filldraw[color=black] (0.5,0.25) circle (2pt);
\filldraw[color=black] (1,0.25) circle (2pt);
\filldraw[color=black] (1.5,0.25) circle (2pt);
\filldraw[color=black] (2,0.25) circle (2pt);
\end{tikzpicture}
\; \; \leftrightarrow \; \;
\begin{tabu}{|[1.3pt]c|c|[1.3pt]c|c|[1.3pt]c|[1.3pt]}
    \tabucline[1.3pt]{-}
     \textcolor{blue}{b} & \textcolor{red}{a} & \textcolor{red}{a} & \textcolor{blue}{b} & \textcolor{red}{a} \\ \hline 
     \textcolor{red}{a} & \textcolor{blue}{b} & \textcolor{blue}{b} & \textcolor{red}{a} & \textcolor{blue}{b}\\ \tabucline[1.3pt]{-}
\end{tabu}
\:  \equiv 1+0+
\begin{array}{|c|}
    \hline
     \textcolor{red}{a} \\ \hline
     \textcolor{blue}{b} \\ \hline
\end{array}
\]
\end{minipage}\vspace{3pt}

\begin{minipage}{1\linewidth}
\[
\begin{tikzpicture}[anchor=base, baseline=-1.5pt]
\draw[green, thick] (0,-0.25) -- (2,-0.25);
\draw[green, thick] (0,0.25) -- (2,0.25);
\draw[green, thick] (0,0.25) -- (0,-0.25) ;
\draw[green, thick] (0.5,0.25) -- (0.5,-0.25) ;
\draw[green, thick] (1,0.25) -- (1,-0.25) ;
\draw[green, thick] (1.5,0.25) -- (1.5,-0.25);
\draw[green, thick] (2,0.25) -- (2,-0.25);
\draw[blue, thick] (1,0.25) -- (1.5,0.25);
\draw[red, thick] (1,-0.25) -- (1.5,-0.25);
\filldraw[color=black] (0,-0.25) circle (2pt);
\filldraw[color=black] (0.5,-0.25) circle (2pt);
\filldraw[color=black] (1,-0.25) circle (2pt);
\filldraw[color=black] (1.5,-0.25) circle (2pt);
\filldraw[color=black] (2,-0.25) circle (2pt);
\filldraw[color=black] (0,0.25) circle (2pt);
\filldraw[color=black] (0.5,0.25) circle (2pt);
\filldraw[color=black] (1,0.25) circle (2pt);
\filldraw[color=black] (1.5,0.25) circle (2pt);
\filldraw[color=black] (2,0.25) circle (2pt);
\end{tikzpicture}
\; \; \leftrightarrow \; \;
\begin{tabu}{|[1.3pt]c|c|[1.3pt]c|c|[1.3pt]c|[1.3pt]}
    \tabucline[1.3pt]{-}
     \textcolor{blue}{b} & \textcolor{red}{a} & \textcolor{blue}{b} & \textcolor{blue}{b} & \textcolor{red}{a} \\ \hline 
     \textcolor{red}{a} & \textcolor{blue}{b} & \textcolor{red}{a} & \textcolor{red}{a} & \textcolor{blue}{b}\\ \tabucline[1.3pt]{-}
\end{tabu}
\:  \equiv 1+3+
\begin{array}{|c|}
    \hline
     \textcolor{red}{a} \\ \hline
     \textcolor{blue}{b} \\ \hline
\end{array}
\]
\end{minipage}\vspace{3pt}

\begin{minipage}{1\linewidth}
\[
\begin{tikzpicture}[anchor=base, baseline=-1.5pt]
\draw[green, thick] (0,-0.25) -- (2,-0.25);
\draw[green, thick] (0,0.25) -- (2,0.25);
\draw[green, thick] (0,0.25) -- (0,-0.25) ;
\draw[green, thick] (0.5,0.25) -- (0.5,-0.25) ;
\draw[green, thick] (1,0.25) -- (1,-0.25) ;
\draw[green, thick] (1.5,0.25) -- (1.5,-0.25);
\draw[green, thick] (2,0.25) -- (2,-0.25);
\draw[blue, thick] (0,-0.25) -- (0,0.25);
\draw[red, thick] (0.5,-0.25) -- (0.5,0.25);
\draw[blue, thick] (1,-0.25) -- (1,0.25);
\draw[red, thick] (1.5,-0.25) -- (1.5,0.25);
\draw[blue, thick] (2,-0.25) -- (2,0.25);
\filldraw[color=black] (0,-0.25) circle (2pt);
\filldraw[color=black] (0.5,-0.25) circle (2pt);
\filldraw[color=black] (1,-0.25) circle (2pt);
\filldraw[color=black] (1.5,-0.25) circle (2pt);
\filldraw[color=black] (2,-0.25) circle (2pt);
\filldraw[color=black] (0,0.25) circle (2pt);
\filldraw[color=black] (0.5,0.25) circle (2pt);
\filldraw[color=black] (1,0.25) circle (2pt);
\filldraw[color=black] (1.5,0.25) circle (2pt);
\filldraw[color=black] (2,0.25) circle (2pt);
\end{tikzpicture}
\; \; \leftrightarrow \; \;
\begin{tabu}{|[1.3pt]c|c|[1.3pt]c|c|[1.3pt]c|[1.3pt]}
    \tabucline[1.3pt]{-}
     \textcolor{blue}{b} & \textcolor{red}{a} & \textcolor{blue}{b} & \textcolor{red}{a} & \textcolor{blue}{b} \\ \hline 
     \textcolor{blue}{b} & \textcolor{red}{a} & \textcolor{blue}{b} & \textcolor{red}{a} & \textcolor{blue}{b}\\ \tabucline[1.3pt]{-}
\end{tabu}
\:  \equiv 4+4+
\begin{array}{|c|}
    \hline
     \textcolor{blue}{b} \\ \hline
     \textcolor{blue}{b} \\ \hline
\end{array}
\]
\end{minipage}
\caption{All 7 (6+1) possible solutions for $2n=10$ players. The first six solutions come from combining blocks $0$, $1$, $2$ and $3$, and the last fixed one comes from block $4$. \textcolor{black}{As per the graphical rules, the nodes at the start and end of the ladder must also have the two lines in green, but here you can see the exception when using block $4$ for odd $n$.}}
\label{fig:casexyx210players}
\end{figure}

On the other hand, if the number of players is not a multiple of 4, i.e. $2n=4k+2$ (equivalently, $n$ is odd), we can argue that we have $k$ blocks of 4 players and a semi-block with 2 players. Then, we will try to relate the number of such solutions to the cases when $2n$ is indeed a multiple of 4 because those last two-player semi-blocks would be fixed. For instance, in the case of 10 players in figure \ref{fig:casexyx210players}, we show all the 7 solutions decomposed as two full blocks and a semi-block. It is worth noticing that in this particular case, the last full block can indeed be $2$ or $3$, as opposed to the situation with full blocks depicted in figure \ref{fig:example12} with 12 players. Therefore, we should in fact re-write the first rule for sticking the blocks as:
\begin{itemize}
    \item \textbf{If n is even}: the solution \emph{cannot start nor end} in blocks $2$ and $3$. 
    \item \textbf{If n is odd}: the solution cannot start with $2$ or $3$ \emph{but it actually can end} in blocks $2$ and $3$.
\end{itemize}

\textcolor{black}{These last two rules can also be translated into graphical rules \mike{as follows:}
For even $n$, the nodes at the end must have the two lines in green, which rules out the possibility of ending  \mike{in block $2$ or $3$}. However, for odd $n$, the last fixed semi-block allows us to get the two green lines from the unmatched lines coming from blocks $2$ and $3$.}

After all the reasoning, we are finally prepared to get all the NE solutions for a generic number of players using $k$ number of combinations of blocks $0$, $1$, $2$ and $3$ \footnote{We ignore block $4$ because there is only one possible combination $4+4+...+4+|b|b|$ when $n$ is odd.}. There is an important remark before starting: firstly, we will build the solutions by combining the blocks \emph{without looking at the starting or ending blocks}, that will result in a recurrence relation, the solution of which needs some initial conditions. The restrictions on the starting/endings blocks will be encoded in such initial conditions. For example, for 6 players ($n=3$ and $k=1$) there are 2 solutions $N(k=1)=2$ and for 10 players ($n=5$ and $k=2$) there are 6 solutions $N(k=2)=6$; whereas for 8 players ($n=4$ and $k=2$) we have $N(k=2)=4$ and for 12 players ($n=6$ and $k=3$) $N(k=3)=10$.

After this clarification, we can proceed to get the solutions. We denote $N_{..j}(k)$, $j\in\{0,1,2,3\}$ as the number of allowed combinations of length $k$ that \emph{end} in block $j$. Then, the total number of chains\footnote{We will be using \emph{combination} and \emph{chain} indistinctly to denote the sequence of blocks that represent one Nash Equilibrium solution.} $N(k)$ is:
\begin{equation}\label{eq:Ntotladcase1}
    N(k)=N_{..0}(k)+N_{..1}(k)+N_{..2}(k)+N_{..3}(k)
\end{equation}
In general, we can always build the chain $N_{..j}(k)$ as a combination of the previous ones by taking into account the rules for sticking the blocks\footnote{In the presented case, the blocks will be connected to the left of the previous one. That's why we use $N_{..j}(k)$ but similarly, the connection can be made to the right, i.e. $N_{j..}(k)$ and the result does not change due to the symmetry of the connection of the blocks.}:
\begin{align}
   & N_{..0}(k)=N_{..0}(k-1)+N_{..1}(k-1)+N_{..3}(k-1) \\
   & N_{..1}(k)=N_{..0}(k-1)+N_{..1}(k-1)+N_{..2}(k-1) \\
   & N_{..2}(k)=N_{..0}(k-1)+N_{..3}(k-1) \\
   & N_{..3}(k)=N_{..1}(k-1)+N_{..2}(k-1) 
\end{align}\label{eq:Nladcase1}
Using that $N_{..2}(k)+N_{..3}(k)=N(k-1)$ and $N_{..0}(k)+N_{..1}(k)=2N(k-1)-N_{..2}(k-1)-N_{..3}(k-1)$ we get to the next recurrence relation:
\begin{equation}\label{eq:recurelation}
N(k)=3N(k-1)-N(k-2)
\end{equation}
Equation (\ref{eq:recurelation}) defines an homogeneous linear recurrence relation with constant coefficients, which can be solved using geometric series, that is, let $N(k)=cr^k$, for some constant $c$. Plugging this ansatz into  (\ref{eq:recurelation}), we arrive to:
\begin{equation}\label{eq:sol1recurelation}
    cr^{k-2}(r^2-3r+1)=0
\end{equation}
which has two non-trivial solutions, $r_1$ and $r_2$:
\begin{equation}\label{eq:r1r2}  
r_1=\frac{1}{2}(3+\sqrt{5})=1+\varphi=\varphi^2\
\; \; \; \; \; \; \; \;
r_2=\frac{1}{2}(3-\sqrt{5})=2-\varphi=\frac{1}{\varphi^2}
\end{equation}
where $\varphi=\frac{1+\sqrt{5}}{2}=1.6180...$ is the golden ratio. In the last equality, we used the well-known identity $\varphi=1+1/\varphi$.

Hence, we can write the general solution for $N(k)$ as:
\begin{equation}\label{eq:genericsol}
N(k)=\alpha r_1^k + \beta r_2^k =\alpha \varphi^{2k}+\beta \varphi^{-2k}
\end{equation}
where $\alpha$ and $\beta$ are set by the initial conditions.

Taking into account the different initial conditions explained above, the relation of $n$ to the number of blocks $k$, and the possible solution using block $4$ for \emph{n even or odd} we can finally write:
\begin{eqnarray}\label{eq:Nfinalxyx2}
    & N(2n)= \frac{2}{\sqrt{5}}\left[\varphi^{n-1} + \varphi^{-(n-1)}\right] & \; \; \; \; \;  n \; even \\
    & N(2n)= \frac{2}{\sqrt{5}}\left[\varphi^{n-1}-\varphi^{-(n-1)}\right] +1 & \; \; \; \; \; n \; odd
\end{eqnarray}
Surprisingly or not, these two equations show that, even though we analysed separately the cases when $n$ is odd or even, the explicit dependence on the golden ratio remains. Note the change of sign in the second term, as well as the extra solution when $n$ is odd. In table \ref{tab:examplescase1} there is a few values for the number of solutions for a given number of players.

\begin{table}
\centering
\caption{The number of solutions for the ladder $N(2n)$ with $2n$ players for the case when the payoff parameters obey $x>y>x/2$.}
\begin{tabular}{|c|c|}
    \hline
    $2n$ & $N_{ladd}(2n)$ \\ \hline
    4 & 2 \\
    6 & 3 \\
    8 & 4 \\
    10 & 7 \\
    12 & 10 \\
    14 & 17 \\
    16 & 26 \\
    18 & 43 \\
    20 & 68 \\ \hline
\end{tabular}
\label{tab:examplescase1}
\end{table}

Clearly, for a large number of players $\varphi^{n-1}$ dominates over $\varphi^{-(n-1)}$, so the approximate solution for both, $n$ even or odd, is:
\begin{eqnarray}\label{eq:Naproxxyx2}
     & N(2n)\approx 0.9\varphi^{n-1}  
\end{eqnarray}

Let's not forget that all of the reasoning in this section (including the next sub-section) is also valid when the parameters obey the relation $2x>y>x$, which is $2(r-p)>q-s>r-p$, since to get to that situation we could perform a \mike{ swap $a\longleftrightarrow b$ } on everything we did, and that would result in the same blocks re-labelled\footnote{Except for block $4$, which would become a different one, reversed, but the same rules as the current one would apply.} with equivalent rules for sticking them. Thus, the result would not change.

\subsection{Circular ladder}\label{sec:closedcase1}
Another interesting set-up is when we close the ladder (known as circular ladder) as shown in figure 4, so there is no ending nor beginning and everyone plays three games. In graph theory, the circular ladder is a polyhedral graph, meaning that is a 3-degree simple planar graph (see [\citet{Gruenbaum2007}] for more on polyhedral graphs).
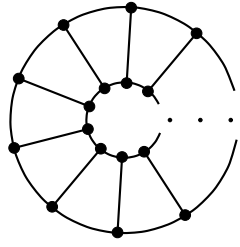
\begin{figure}
\centering
\begin{tikzpicture}
\draw[black, thick, domain=25:340] plot[smooth] ({0.5*cos(\x)}, {0.5*sin(\x)});
\draw[black, thick, domain=15:350] plot[smooth] ({1.5*cos(\x)}, {1.5*sin(\x)});
\filldraw[color=black, thin] (0.6,0) circle (0.65pt);
\filldraw[color=black, thin] (1.0,0) circle (0.65pt);
\filldraw[color=black, thin] (1.4,0) circle (0.65pt);
\filldraw[color=black] (0.31,0.38) circle (2pt);
\filldraw[color=black] (0.03,0.49) circle (2pt);
\filldraw[color=black] (-0.26,0.42) circle (2pt);
\filldraw[color=black] (-0.46,0.18) circle (2pt);
\filldraw[color=black] (-0.48,-0.12) circle (2pt);
\filldraw[color=black] (-0.31,-0.38) circle (2pt);
\filldraw[color=black] (-0.03,-0.49) circle (2pt);
\filldraw[color=black] (0.26,-0.42) circle (2pt);
\filldraw[color=black] (0.95,1.15) circle (2pt);
\filldraw[color=black] (0.09,1.49) circle (2pt);
\filldraw[color=black] (-0.80,1.26) circle (2pt);
\filldraw[color=black] (-1.39,0.55) circle (2pt);
\filldraw[color=black] (-1.45,-0.37) circle (2pt);
\filldraw[color=black] (-0.95,-1.15) circle (2pt);
\filldraw[color=black] (-0.09,-1.49) circle (2pt);
\filldraw[color=black] (0.80,-1.26) circle (2pt);
\draw[black, thick] (0.31,0.38) -- (0.95,1.15);
\draw[black, thick] (0.03,0.49) -- (0.09,1.49);
\draw[black, thick] (-0.26,0.42) -- (-0.80,1.26);
\draw[black, thick] (-0.46,0.18) -- (-1.39,0.55);
\draw[black, thick] (-0.48,-0.12) -- (-1.45,-0.37);
\draw[black, thick] (-0.31,-0.38) -- (-0.95,-1.15);
\draw[black, thick] (-0.03,-0.49) -- (-0.09,-1.49);
\draw[black, thick] (0.26,-0.42) -- (0.80,-1.26);
\end{tikzpicture}
\caption{Circular ladder graph. In this completely symmetric situation, everyone plays three games.}
\label{fig:closeladd}
\end{figure}

We could also ask how many NE there are in the new situation. Fortunately, the previous analysis with the (open) ladder helps us, since the circular ladder is just a ladder in which the last two players are connected to the first two. For simplicity, we will only consider the case when the number of players $2n$ is a multiple of 4 (\textbf{$n$ is even}), but the procedure is extendable to the case when $n$ is odd by using the semi-block argument and as we saw previously, the explicit dependence for a large number of players might not change. 

Taking into account that we will only be working with full blocks and that no block is left unconnected, we can make a distinction between the first four blocks ($0$, $1$, $2$ and $3$) and the last one ($4$). Let's recall that block $4$ can only be connected with itself, which now, for even $n$ and the circular ladder, adds only one fixed combination ($44444...4$). Moreover, swapping $a$ and $b$ in block $4$ gives an additional extra block only allowed for this particular case, that, again, can only be connected with itself, adding a second fixed possible solution. Therefore, to the total number of allowed combinations for the circular ladder using blocks $0$, $1$, $2$ and $3$, we need to add two solutions coming from block $4$ and its reversed counterpart.

Having specified that subtlety, we can move on to the analysis using the first 4 blocks. The first thing to notice is that, as opposed to the ladder, there are no restrictions on the ending or starting blocks, but not every ending block can be attached to the starting one. For example, in the current case, the combination starting with $0$ and ending in $2$ (i.e. $0...2$) is not a valid solution since those two blocks cannot be connected, whereas the reversed ones can (i.e. $2...0$). \textcolor{black}{In graphical terms, the rule of every node having at least two green lines applies, and connecting the ``last nodes" to the starting ones gives more freedom to get (at least) those two green lines}. This makes the analysis more tricky, but it still can be done; we just need to count all the possible combinations whose ending block can be attached to the first one. From the previous rules, there are 10 allowed combinations out of the 16 possible, which are:
\begin{itemize}
    \item[•] $0...0$ ; $0...1$ ; $0...3$
    \item[•] $1...0$ ; $1...1$ ; $1...2$
    \item[•] $2...0$ ; $2...3$
    \item[•] $3...1$ ; $3...2$
\end{itemize}
In order to do the counting, we would need to perform a similar analysis for building the ladder as in equations (\ref{eq:Ntotladcase1}) and (8)-(11), but now from right to left. That is, now, we would consider $N_{j...}(k)$, $j \in \{0,1,2,3\}$ as the number of allowed combinations of length $k$ that \emph{start} with $j$. We would arrive to the same recurrence relation as in (\ref{eq:recurelation}) but now, since we are not restricting the starting or ending blocks, the initial conditions change to $N(k=1)=4$ and $N(k=2)=10$ \footnote{Those conditions give the same list of numbers as in equation (\ref{eq:Nfinalxyx2}) for even $n$ but shifted ahead by one position.}. So, the new solution is:
\begin{equation}\label{eq:closed1case1}
N(k)=\frac{2}{\sqrt{5}}\left[\varphi^{2k+1}+\varphi^{-2k-1}\right]
\end{equation}

Equation (\ref{eq:closed1case1}) gives the number of combinations using blocks $0$, $1$, $2$ and $3$ \emph{without any restriction} on the starting or ending blocks. We will use this solution later on, because now we are primarily interested in getting the number of combinations that start in block $j$ and end in block $m$, that is $N_{j..m}(k)$. To do so, we can analyse the chain when building it to the left with $N_{..m}(k)$ and to the right with $N_{j..}(k)$.

From the relations in equations (8)-(11) we can write the number of combinations starting or ending in $0$ or $1$ as a function of the number that start or end in $2$ and/or $3$:
\begin{align}\label{eq:relationscase1}
    & N_{0..}(k)=N_{2..}(k)+N_{3..}(k)-N_{3..}(k-1) \\
    & N_{1..}(k)=N_{2..}(k)+N_{3..}(k)-N_{2..}(k-1) \\
    & N_{..0}(k)=N_{..2}(k)+N_{..3}(k)-N_{..2}(k-1) \\
    & N_{..1}(k)=N_{..2}(k)+N_{..3}(k)-N_{..3}(k-1) 
\end{align}

All the relations above tell us that we can focus only on getting $N_{2..}(k)$, $N_{3..}(k)$, $N_{..2}(k)$ and $N_{..3}(k)$ because all the others can be obtained as a function of those.
There is also a relation between $N_{2..}(k)$ and $N_{3..}(k)$ \footnote{Same relation holds for $N_{..2}(k)$ and $N_{..3}(k)$.}:
\begin{equation}
    N_{2..}(k)+N_{3..}(k)=N(k-1) 
\label{eq:closecase12and2}
\end{equation}
which, by symmetry\footnote{The recurrence relation depending only on $N_{2..}(k)$ and $N_{3..}(k)$ would be identical, so the total number $N(k-1)$ must be equally split between $N_{2..}(k)$ and $N_{3..}(k)$.}, leads to $N_{2..}(k)=N_{3..}(k)=N(k-1)/2$ and similarly, $N_{..2}(k)=N_{..3}(k)=N(k-1)/2$. We also have $N_{..2}(k)=N_{2..}(k)$ and $N_{..3}(k)=N_{3..}(k)$. Expanding $N_{2..}(k)$ and $N_{..2}(k)$ and using the relations in (20)-(23):
\begin{align}
    & N_{2..}(k)=3N_{2..2}(k)-N_{2..2}(k-1)+3N_{2..3}(k)-N_{2..3}(k-1) \\
    & N_{..2}(k)=3N_{2..2}(k)-N_{2..2}(k-1)+3N_{3..2}(k)-N_{3..2}(k-1) 
\end{align}
\label{eq:closedrel2..}

These last two equations lead to $N_{2..3}(k)=N_{3..2}(k)$. Identical relations hold for $N_{3..}(k)$ and $N_{..3}(k)$, so we can write $N_{2..2}(k)=N_{3..3}(k)$.

Moreover, in equations (25) and (26) the first two terms for $N_{2..2}(k)$ and the last two with $N_{2..3}(k)$ have the form of the recurrence relation in equation (\ref{eq:recurelation}), so rearranging\footnote{The second equation with $N_{..2}(k)$ would give the same information, so we omit it.}:
\begin{equation}
N_{2..}(k)=N_{2..2}(k+1)+N_{2..3}(k+1)
\label{eq:closedrel2223}
\end{equation}
which gives us $N_{2..2}(k)=N_{2..3}(k)=N_{2..}(k-1)/2$ and knowing $N_{2..}(k)=N(k-1)/2$, we are finally ready to write:
\begin{equation}
N_{2..2}(k)=N_{2..3}(k)=\frac{N(k-2)}{4}
\label{eq:closedrel2N}
\end{equation}

To summarise, after all the reasoning and relations, we arrived to the conclusion that we only need to know $N_{2..2}(k)$ because from that we can obtain the number for any other combination $N_{j..m}(k)$ by using the relations in (20)-(23). After adding the 10 allowed ones, plus the two solutions from blocks $4$ and the reversed one and using the recurrence relation in (\ref{eq:recurelation}) we can get the total number of solutions for the circular ladder (see the section \ref{sec:appendix1} of the appendix for the details of the calculation and a subtlety regarding the final $+2$):
\begin{equation}\label{eq:closeladsolutioncase1}
    N_{circ}(k)=2N(k-1)-\frac{N(k-2)}{2}+2
\end{equation}
where $N(k)$ is given in (\ref{eq:closed1case1}). The final expression for $N_{circ}(k)$ is then:
\begin{equation}\label{eq:finalcloseladsolutioncase1}
    N_{circ}(k)=\varphi^{2k}+\varphi^{-2k}+2
\end{equation} 
Obviously, for large $k$ and using $k=n/2$, then  (\ref{eq:finalcloseladsolutioncase1}) becomes\footnote{For $k>5$, that is for more than 20 players, the relative error of the approximate $N(2n)$ is below $1\%$}:
\begin{equation}\label{eq:aproxclosesolcase1}
    N_{circ}(2n)\approx \varphi^{n} \; \; \; \; \; \; n \; \; even
\end{equation}

\begin{figure}
    \centering
    \includegraphics[scale=0.65]{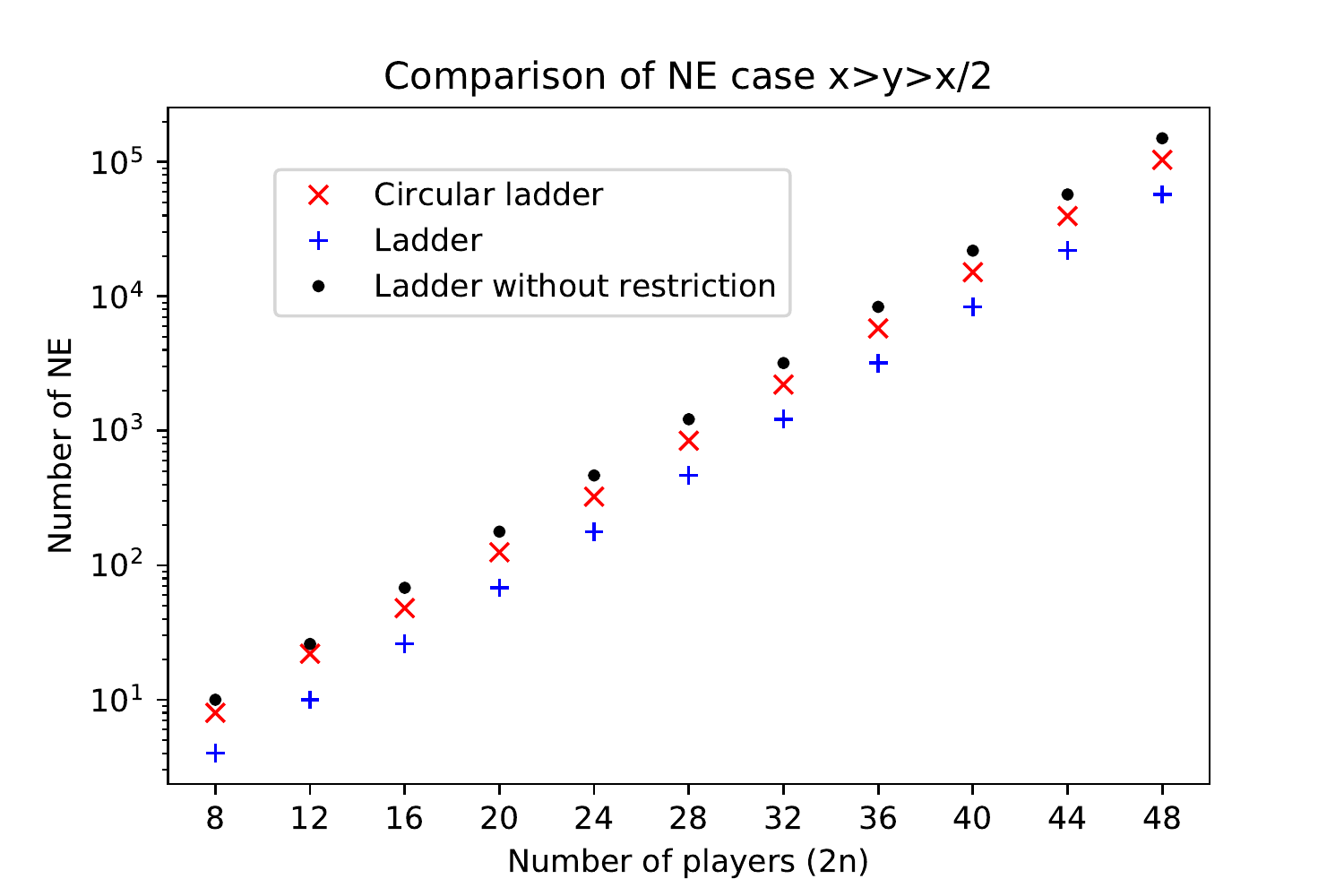}
    \caption{The number of Nash Equilibria as a function of the number of players $2n$ (assumed to be a multiple of four, i.e. even $n$) for the ladder, the circular ladder and the ladder without any ending/beginning restriction when the parameters obey $x>y>x/2$ (translated into the original parameters, $(r-p)>(q-s)>(r-p)/2$). Notice the logarithmic scale used for the number of NE.}
    \label{fig:plotcase1}
\end{figure}

It is interesting to compare this result for the circular ladder $N_{circ}(2n)$ with the result in (\ref{eq:Naproxxyx2}) for the ladder $N_{ladd}(2n)$ and in  (\ref{eq:closed1case1}) for the ladder without beginning/ending restriction $N_{lnr}(2n)$, as shown in figure \ref{fig:plotcase1}. The number of NE grows exponentially for all of them as $\varphi^n\approx (1.62)^n$ (for even $n$), but with a different scaling factor. The relations are: $N_{circ}(2n)\approx 1.80N_{ladd}(2n)$, $N_{circ}(2n)\approx 0.69N_{lnr}(2n)$, and $N_{ladd}(2n)\approx 0.38N_{lnr}(2n)$

As we argued in the previous sub-section, all the reasoning applies as well to the situation when $2x>y>x$, that is $2(r-p)>q-s>r-p$.

\section{Case when $x/2>y$}\label{sec:casex2y}
In this section we will be performing the same analysis as in the previous section, but now the parameters obey the relation $x/2>y$, which means $r-p>2(q-s)$.

The rule for the best strategy in this case is: \emph{play $a$ only when everyone else plays $b$, and play $b$ otherwise}\footnote{Equivalently, if at least one player chooses $a$, then play $b$.}. This case presents a richer number of combinations, so there are actually six different building blocks to be combined with each other:

\[
0 \equiv 
\begin{array}{|c|c|}
    \hline
     \textcolor{red}{a} & \textcolor{blue}{b}  \\ \hline
     \textcolor{blue}{b} & \textcolor{red}{a} \\ \hline
\end{array}
\; \; \; \;  
1 \equiv 
\begin{array}{|c|c|}
    \hline
     \textcolor{red}{a} & \textcolor{blue}{b}    \\ \hline
     \textcolor{blue}{b} & \textcolor{blue}{b} \\ \hline
\end{array}
\; \; \; \;  
2 \equiv 
\begin{array}{|c|c|}
    \hline
     \textcolor{blue}{b} & \textcolor{red}{a}  \\ \hline
     \textcolor{red}{a} & \textcolor{blue}{b} \\ \hline
\end{array}
\; \; \;  \; 
3 \equiv 
\begin{array}{|c|c|}
    \hline
     \textcolor{blue}{b} & \textcolor{blue}{b}  \\ \hline
     \textcolor{red}{a} & \textcolor{blue}{b} \\ \hline
\end{array}
\; \; \; \;  \; 
4 \equiv 
\begin{array}{|c|c|}
    \hline
     \textcolor{blue}{b} & \textcolor{red}{a}  \\ \hline
     \textcolor{blue}{b} & \textcolor{blue}{b} \\ \hline
\end{array}
\; \; \; \; 
5 \equiv 
\begin{array}{|c|c|}
    \hline
     \textcolor{blue}{b} & \textcolor{blue}{b}  \\ \hline
     \textcolor{blue}{b} & \textcolor{red}{a} \\ \hline
\end{array}
\]

We will provide the rules for the linking between each block when the connection is done to the left and to the right of the block, which are in table \ref{tab:rules}. For instance, after $0$ (to the right of it), there can only be attached $0$, $1$ or $4$ and before $0$ (to the left), there can only be $0$, $3$ or $5$.

\begin{table}[]
    \centering
    \caption{Rules for connecting the allowed blocks to the left (first column) and to the right (last column) of a given block (in bold).}
    \begin{tabular}{|ccccc|}
        \hline
       $0, 3, 5$  & $\rightarrow$ & $\textbf{0}$ & $\rightarrow$ & $0, 1, 4$  \\ 
        $0, 3, 5$ & $\rightarrow$ & $\textbf{1}$ & $\rightarrow$ &  $3, 2$  \\ 
        $1, 2, 4$ & $\rightarrow$ & $\textbf{2}$ & $\rightarrow$ &  $2, 3, 5$  \\ 
        $1, 2, 4$ & $\rightarrow$ & $\textbf{3}$ & $\rightarrow$ &  $0, 1$  \\ 
        $0, 5$ & $\rightarrow$ & $\textbf{4}$ & $\rightarrow$ &  $2, 3, 5$  \\ 
        $2, 4$ & $\rightarrow$ & $\textbf{5}$ & $\rightarrow$ &  $0, 1, 4$  \\ \hline
    \end{tabular}
    \label{tab:rules}
\end{table}

From these rules, it is worth noticing the symmetry for blocks $0$ and $2$: there are three possible blocks attached to the left and to the right. That leads to the expectation of having the same number of chains that start and end with $0$, and similarly, for block $2$. Nonetheless, that is not the case for the other blocks. For example for block $1$, there are only two possibilities after $1$, but three before it. This asymmetry leads to the conclusion that there won't be the same number of blocks starting and ending in $1$. Same argument applies to blocks $3$, $4$ and $5$.

\textcolor{black}{As before, the \mike{graphical rules} are more revealing. Again, a red line around a node/player means choosing $a$, a blue one means $b$; and connecting both lines translates into a green one. This time, the graphical rule changes slightly; while in the previous section each node must have at least two green lines, in this case, each node \mike{must only have at 
least one} green line. Such \mike{a} difference in the graphical rule \mike{{\it a priori} tells} us that in the present case there will be more NE than in the previous one.} 
\subsection{Ladder}\label{sec:laddercase2}
Having the rules and blocks, we are able to get all the solutions. For simplicity with the forthcoming analysis, we will assume that the number of players is a multiple of 4 (\textbf{n is even}). As previously, there are certain restrictions on the starting and ending blocks: the solution cannot start with $4$ nor $5$ and cannot end in $1$ nor $3$. \textcolor{black}{Otherwise, \mike{ either the starting or ending node, respectively, would not have the minimum of one green line required.}}
The initial conditions in this case are the same as in the previous case for \emph{odd} $n$: the number of solutions for only one block is two $N(k=1)=2$ (chains ``$0$", and ``$2$"), and for two blocks, there are only six $N(k=2)=6$ (combinations ``$00$", ``$04$", ``$12$", ``$22$", ``$25$" and ``$30$").

As before, let's denote by $N_{..j}(k)$, $j\in\{0,1,2,3,4,5\}$ the number of allowed combinations of length $k$ that end in $j$. Then, the total number of chains $N(k)$ is:
\begin{equation}\label{eq:Ntotladcase2}
    N(k)=N_{..0}(k)+N_{..1}(k)+N_{..2}(k)+N_{..3}(k)+N_{..4}(k)+N_{..5}(k)
\end{equation}

We can always build the chain $N_{..j}(k)$ as a combination of the previous ones by taking into account the rules mentioned above:
\begin{align}
   & N_{..0}(k)=N_{..0}(k-1)+N_{..3}(k-1)+N_{..5}(k-1) \\
   & N_{..1}(k)=N_{..0}(k-1)+N_{..3}(k-1)+N_{..5}(k-1) \\
   & N_{..2}(k)=N_{..1}(k-1)+N_{..2}(k-1)+N_{..4}(k-1) \\
   & N_{..3}(k)=N_{..1}(k-1)+N_{..2}(k-1)+N_{..4}(k-1) \\
   & N_{..4}(k)=N_{..0}(k-1)+N_{..5}(k-1) \\
   & N_{..5}(k)=N_{..2}(k-1)+N_{..4}(k-1)
\end{align}\label{eq:Nladcase2}
Using that $N(k-1)=N_{..0}(k)+N_{..2}(k)=N_{..1}(k)+N_{..3}(k)$ we get to the exact same recurrence relation as in (\ref{eq:recurelation}), with the same values for $\alpha$ and $\beta$ in the case when $n$ was odd. Hence, keeping the dependence in $k$:
\begin{equation}\label{eq:solrecurelat}
N(k)=\frac{2}{\sqrt{5}}(\varphi^{2k}-\varphi^{-2k})
\end{equation}

As earlier, for a large number of players and using the relation between $k$ and the number of players $2n=4k$ then:
\begin{equation}
  N(2n)\approx 0.9\varphi^n \; \; \; \; \; \; n \; \; even
\end{equation}

As we did in the previous section, we could try to get the number of NE in the case when the number of players is not a multiple of 4, but it's still a tedious task, and, instead, we will focus on the analysis for the circular ladder and the comparison with the results on the previous case.

\subsection{Circular ladder}\label{sec:closedcase2}
In this section, we will perform the same analysis as with the circular ladder in the previous section. Now, from the rules in table \ref{tab:rules}, we can extract the 16 allowed combinations, \textcolor{black}{in which every starting and ending node has at least one green line. Using blocks, such combinations are}:
\begin{itemize}
    \item[•] $0...0$ ; $0...3$ ; $0...5$
    \item[•] $1...0$ ; $1...3$ ; $1...5$
    \item[•] $2...1$ ; $2...2$ ; $2...4$
    \item[•] $3...1$ ; $3...2$ ; $3...4$
    \item[•] $4...0$ ; $4...5$
    \item[•] $5...2$ ; $5...4$
\end{itemize}

The initial conditions for the solution without any starting/ending restriction are $N(k=1)=6$ and $N(k=2)=16$ \footnote{As earlier, those conditions give the same list of numbers as in equation (\ref{eq:solrecurelat}) but shifted ahead by one position.}. The new solution is:
\begin{equation}\label{eq:closesolrecurelat}
N(k)=\frac{2}{\sqrt{5}}\left[\varphi^{2k+2}-\varphi^{-2k-2}\right]
\end{equation}

From the relations in equations (33)-(38) we can write $N_{..1}(k)=N_{..0}(k)$ and $N_{..3}(k)=N_{..2}(k)$, and similarly, for the left-to-right analysis we get $N_{4..}(k)=N_{2..}(k)$ and $N_{5..}(k)=N_{0..}(k)$. Those relations imply that there will be the same number of chains ending in  block $1$ (block $3$) and ending in block $0$ (block $2$), and the same number starting with $4$ and $0$ ($5$ and $2$).
For the remaining ones, we get the relations:
\begin{align}
    & N_{1..}(k)=N_{2..}(k)-N_{0..}(k-1) \\
    & N_{3..}(k)=N_{0..}(k)-N_{2..}(k-1) \\
    & N_{..4}(k)=N_{..0}(k)-N_{..2}(k-1)  \\
    & N_{..5}(k)=N_{..2}(k)-N_{..0}(k-1)  
\end{align} \label{eq:closedcombin}

We will omit all the tedious math because the procedure has been explained in section \ref{sec:closedcase1}. We arrive to analogous equations such as $N_{2..0}(k)=N_{0..2}(k)$ ; $N_{0..0}(k)=N_{2..2}(k)$ and the key relation:
\begin{equation}\label{eq:closedrel0N}
N_{0..0}(k)=N_{0..2}(k)=\frac{N(k-2)}{4}
\end{equation}

Once we have this last equality, we can use $N_{0..0}(k)$ to get any other combination $N_{j..m}(k)$ by using the relations in (42)-(45). After adding the 16 allowed combinations and using the recurrence relation in (\ref{eq:recurelation}) we can get the total number of solutions for the circular ladder (see the section \ref{sec:appendix2} of the appendix for the details of the calculation):
\begin{equation}
    N_{circ}(k)=3N(k-2)-\frac{N(k-4)}{2}
\label{eq:closeladsolution}
\end{equation}
where $N(k)$ is given in equation (\ref{eq:closesolrecurelat}). Finally, the expression for $N_{circ}(k)$ is:
\begin{equation}
    N_{circ}(k)=\varphi^{2k}+\varphi^{-2k}
\label{eq:finalcloseladsolution}
\end{equation}
This result might be surprising, because, except for the $+2$, it's the same solution as the one we found for the circular ladder in the previous section in (\ref{eq:finalcloseladsolutioncase1}). Clearly, for a large number of players, both solutions for the circular ladder are identical:
\begin{equation}\label{eq:aproxclosesol}
    N_{circ}(2n)\approx \varphi^n \; \; \; \; \; \; n \; \; even
\end{equation}

\begin{figure}
    \centering
    \includegraphics[scale=0.65]{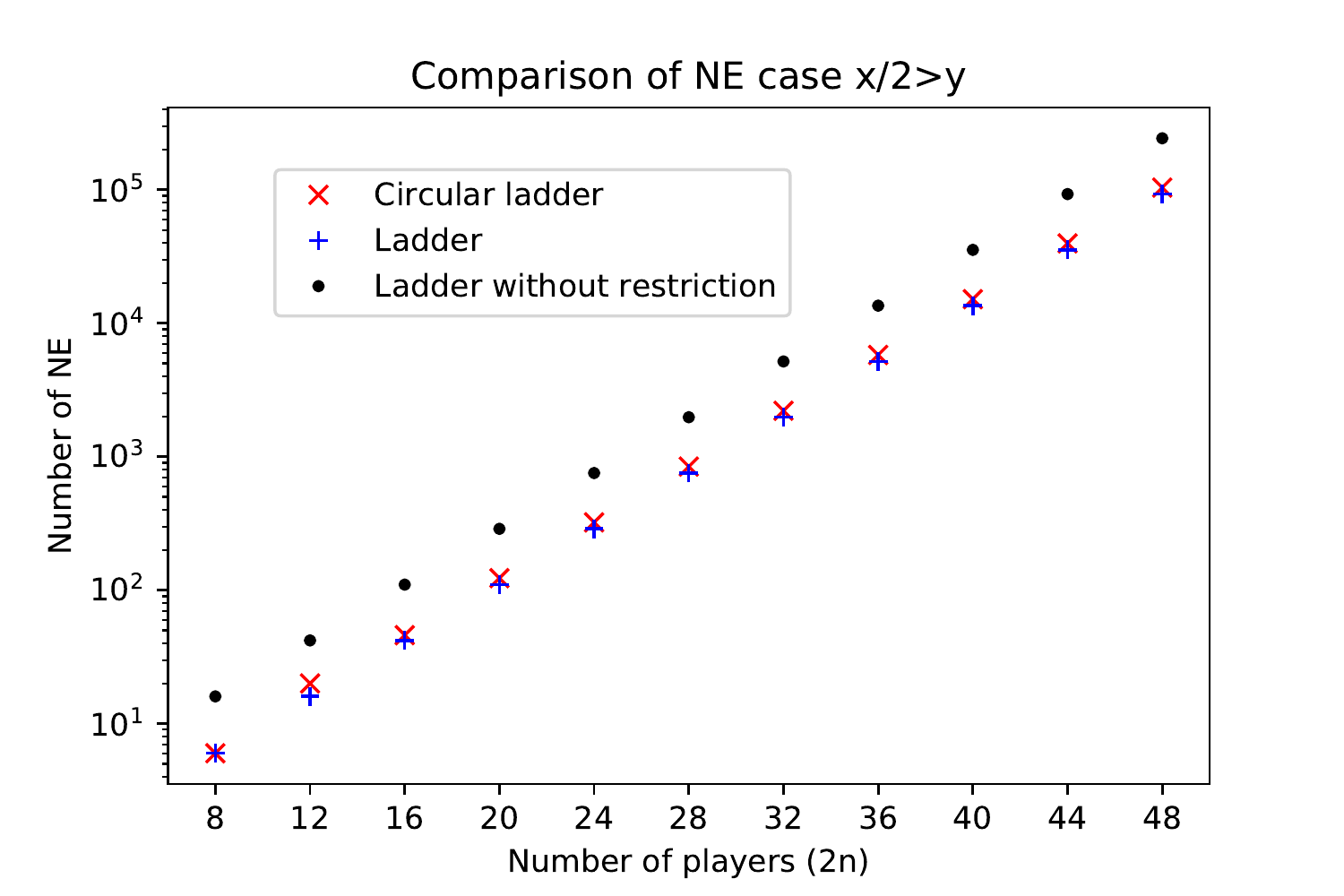}
    \caption{The number of Nash Equilibria (in logarithmic scale) as a function of the number of players $2n$ (even $n$) for the ladder, the circular ladder and the ladder without any ending/beginning restriction when the parameters obey $x/2>y$ (in the original parameters, $(r-p)/2>(q-s)$).}
    \label{fig:plotcase2}
\end{figure}

Figure \ref{fig:plotcase2} shows the comparison between the ladder, the circular ladder and the ladder without beginning/ending restrictions. As before, the ladder without restrictions presents the largest number of NE: $N_{ladd}(2n)\approx 0.38N_{lnr}(2n)$ and $N_{circ}(2n)\approx 0.43N_{lnr}(2n)$. This time, there is not such a big difference between the ladder and the circular ladder $N_{circ}(2n)\approx 1.12N_{ladd}(2n)$, only $12\%$ percent; in contrast to the previous case, in which the ratio reached $80\%$.

Again, all the reasoning applies as well to the situation when $y>2x$, that is $q-s>2(r-p)$, because swapping $a\longleftrightarrow b$ would result in different blocks, but with equivalent rules for attaching them, leaving our result unchanged.

\section{Discussion}\label{sec:discussion}
\begin{figure}
    \centering
    \includegraphics[scale=0.65]{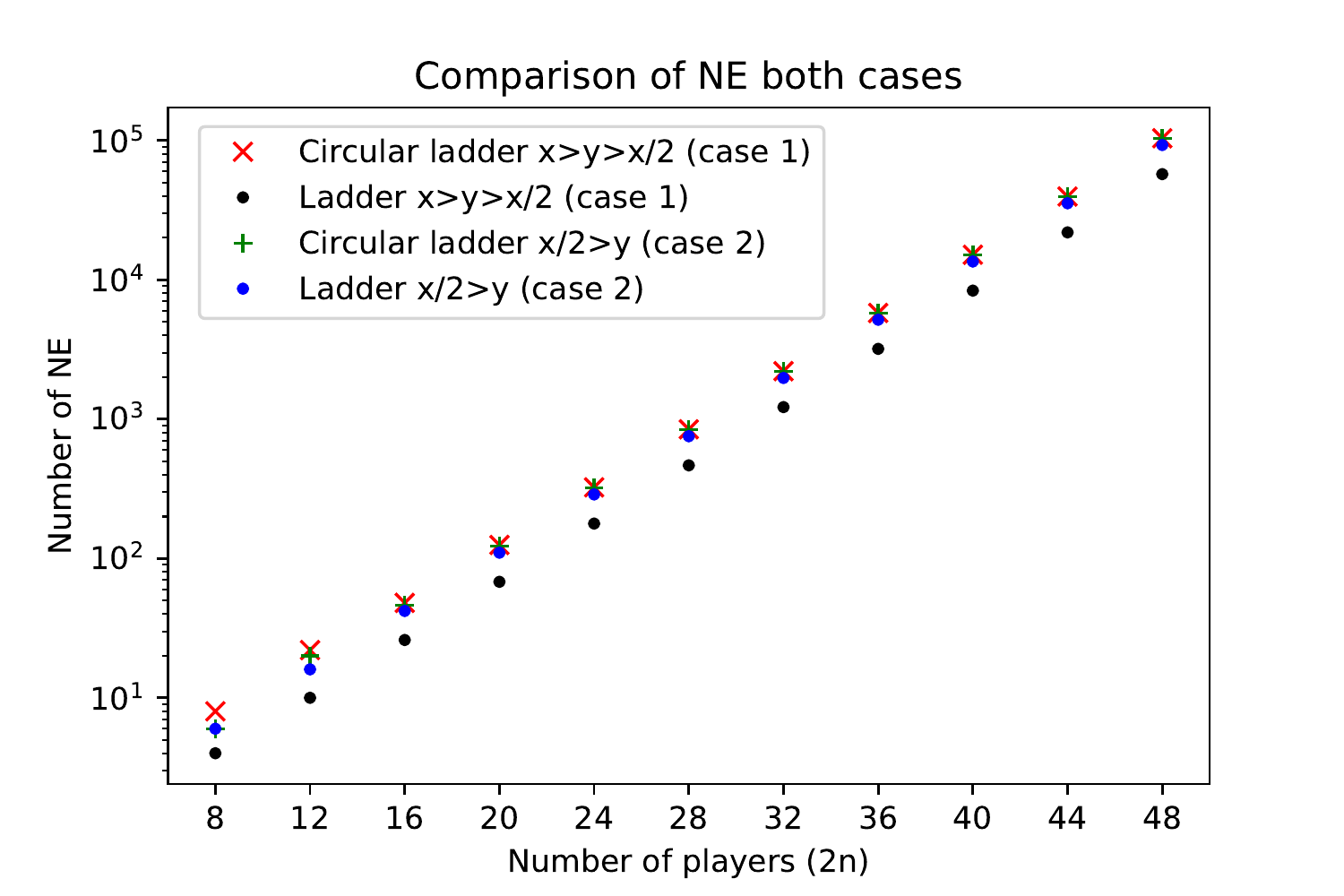}
    \caption{Comparison of the results for the number of Nash Equilibria as a function of the number of players $2n$ for the ladder and the circular ladder for the two cases discussed. For a large number of players, the red cross marker corresponding to the circular ladder for case 1 overlaps the green plus marker for the circular ladder for case 2 (all along, the difference between both is only two solutions).}
    \label{fig:plotcomparison}
\end{figure}

In figure \ref{fig:plotcomparison}, we compare the results of our two representative cases for the ladder and the circular ladder\footnote{For the discussion, we omit the ladder without restrictions because it's just an intermediate result with no further interest.}. A quick look reveals that which relation the payoff parameters obey indeed affects the results for the ladder: $N_{ladd}^{case 2}(2n)\approx 1.62N_{ladd}^{case 1}(2n)$ (black and blue dots in figure \ref{fig:plotcomparison}). Hence, case 2 has about $62\%$ more solutions than case 1 for the ladder, \textcolor{black}{as could have been inferred by the difference in the graphical rule}. \mike{This} distinction of cases is not important for the circular ladder, since $N_{circ}^{case 1}(2n)\approx N_{circ}^{case 2}(2n)$ (the overlapping red cross marker and green plus marker) \textcolor{black}{even when the graphical rules were still different. Such \mike{an} outcome might suggest} that, for this particular case of the 3-degree graph, the number of NE is directly related to the players' settings (i.e. the graph itself) and the definition of the NE for the two-player game.

Summarising, we found that, despite the value of the payoff parameters, but still assuming the two NE for the two-player game, and up to a multiplying constant, the number of solutions for the ladder and the circular ladder is exponential in (half) the total number of players $n$, as $\varphi^n \approx (1.62)^n$ (for even $n$). The scaling factor is different for the ladder and the circular ladder when the payoff parameters obey one particular relation (as discussed, if case 1 or case 2). Moreover, such factor is affected by the value of the payoff parameters for the ladder setting (the number of NE for the ladder changes when considering case 1 or case 2), while such payoff parameters do not change the number of NE for the circular ladder. It is also remarkable the ubiquity of the golden ratio $\varphi$ through the whole analysis for both graphs.

\textcolor{black}{To our knowledge, no example of a graphical game with a large number of players has been approached analytically \mike{as we have done}
with this particular example. It's also worth mentioning that throughout our procedure an algorithmic approach can be seen in the sense that our 4-player blocks represent a meta-player, which can (or not) be matched with \mike{other} meta-players in a purely combinatoric way. \mike{ The associated rules } can also be translated into certain colouring rules on the graph. Therefore, our procedure could \mike{probably} be extended to other graphical games with different (highly regular and possibly complete) underlying graphs.}

\section{Further study}\label{sec:further}
\textcolor{black}{It would be interesting to see whether we could reformulate the problem
\mike{in the following way:} \mike{instead} of having $2n$ players with two possible strategies and $3n-2$ copies \mike{of} the two-player game (in the case of the ladder, for the circular there are $3n$), \mike{consider} for example two players playing only one equivalent game but with many more strategies available, which has been broadly studied over the years (for example in [\citet{MCLENNAN2005264}]).}

\textcolor{black}{Another further step would be studying the generic situation of playing three completely different anti-coordination games with each player, which would require taking into account the orientation/direction of the graph, as well as trying to obtain the Pareto Optimal solution/s among all the NE and compare it to the other NE.}

\textcolor{black}{In the type of game studied, it might be worth doing an equivalent analysis with a coordination game, with one or two NE with symmetric and asymmetric payoff parameters, or even a \mike{mix of } coordination and anti-coordination games (e.g. play a coordination game with the player in front and the same anti-coordination game with the players at the sides).}

\textcolor{black}{In changing the type of graph, it \mike{could} also be \mike{illuminating} to study the same anti-coordination game for a} degree 4 graph, for instance, a graph in which there is a regular outer polygon, and the same polygon inside rotated, and each vertex is connected to the closest 4 vertices (one example could be the pentacle). In doing do, we would be interested in knowing if the number of NE solutions can be generalised for a specific or generic n-degree graph.

\section*{Acknowledgements}
\mike{We would like to thank the editor and a referee for constructive comments on an original version of this paper.}

\newpage


\appendix
\numberwithin{equation}{section}

\section{Circular ladder when $x>y>x/2$}\label{sec:appendix1}
The total number of allowed combinations for the circular ladder are:
\begin{equation}
\begin{aligned}
    N_{circ}(k) & =  N_{0..0}(k)+N_{0..1}(k)+N_{0..3}(k)+N_{1..0}(k)+N_{1..1}(k)+\\ 
    & + N_{1..2}(k)+N_{2..0}(k)+N_{2..3}(k)+N_{3..1}(k)+N_{3..2}(k)
\end{aligned}\label{eq:Nclosedtotalcase1}
\end{equation}
We can simplify (\ref{eq:Nclosedtotalcase1}) by using the fact that the allowed combinations are the total minus the non-allowed ones. Then:
\begin{equation}
\begin{aligned}
N_{circ}(k)=& N(k)-N_{0..2}(k)-N_{1..3}(k)-N_{2..1}(k)-\\
 &-N_{2..2}(k)-N_{3..0}(k)-N_{3..3}(k)
\end{aligned}\label{eq:Nclosedtotalcase102}
\end{equation}

We can rewrite each term as a function of $N_{2..2}(k)$ and $N_{2..3}(k)$ by using the relations in (20)-(23) and $N_{2..3}(k)=N_{3..2}(k)$:
\begin{align}
    & N_{3..3}(k)=N_{2..2}(k) \\
    & N_{0..2}(k)=N_{2..2}(k)+N_{2..3}(k)-N_{2..3}(k-1) \\
    & N_{3..0}(k)=N_{2..2}(k)+N_{2..3}(k)-N_{2..3}(k-1) \\
    & N_{1..3}(k)=N_{2..2}(k)+N_{2..3}(k)-N_{2..3}(k-1) \\
    & N_{2..1}(k)=N_{2..2}(k)+N_{2..3}(k)-N_{2..3}(k-1) 
\end{align}\label{eq:Nasfunction2223}

Adding up all of those, and rearranging, (\ref{eq:Nclosedtotalcase102}) becomes:
\begin{equation}\label{eq:Nclosedtotalcase12almost}
\begin{aligned}
    N_{circ}(k) &= N(k)-4\left[N_{2..2}(k)+N_{2..3}(k)\right]-2\left[N_{2..2}(k)-2N_{2..3}(k-1)\right] \\
    &= N(k)-2N(k-2)-2\left[N_{2..2}(k)-2N_{2..3}(k-1)\right]
\end{aligned}
\end{equation}
in second equality we used:
\begin{equation}\label{eq:N22N23}
N_{2..2}(k)+N_{2..3}(k)=N_{2..}(k-1)=\frac{N(k-2)}{2}   
\end{equation}

Up to this point, a remark is worth mentioning: we are counting the number of combinations that end and start with some given blocks, so $N_{j..m}(k)$ must be an integer. From (\ref{eq:N22N23}) we could write $N_{2..2}(k)=N_{2..3}(k)=N(k-2)/4$, but for certain $k$ it would lead to non-integer numbers of combinations. For example, for the current case, the counting of the combinations gives that for $k=5$, $N_{2..2}(5)$ is one less than $N_{2..3}(5)$ [$N_{2..2}(5)=6$ and $N_{2..3}(5)=7$]; for $k=6$ they are both equal [$N_{2..2}(6)=N_{2..3}(6)=17$] and for $k=7$, now $N_{2..2}(7)$ has one more combination than $N_{2..3}(7)$ [$N_{2..2}(7)=45$ and $N_{2..3}(7)=44$] \footnote{For the next $k$, the cycle would start again, being $N_{2..2}(8)$ one below $N_{2..3}(8)$, then for $k=9$ equal, and finally, for $k=10$ one above.}. The specific multiplying terms for $N_{2..2}(k)$ and for $N_{2..3}(k)$ lead to the fact that $2/3$ of the time we would be one above the actual counting, and $1/3$ two below\footnote{This unequal ratio comes from $N(k)$ given in (\ref{eq:closed1case1}), because $2/3$ of the time, $N(k)/4$ would be an odd number over 2, and the rest $1/3$ would be an integer.}.

Nonetheless, after having clarified this, we can plug $N_{2..2}(k)=N_{2..3}(k)=N(k-2)/4$ into (\ref{eq:Nclosedtotalcase12almost}) to get to:
\begin{equation}
\begin{aligned}
    N_{circ}(k)&= N(k)-\frac{5}{2}N(k-2)+N(k-3) \\
    &= N(k)-N(k-1)+\frac{N(k-2)}{2} \\
    &= 2N(k-1)-\frac{N(k-2)}{2} 
\end{aligned}\label{eq:Nclosedtotalcase12}
\end{equation}
In the last two lines we used the recurrence relation $N(k)=3N(k-1)-N(k-2)$.
Remember that this solution is for the combination of blocks $0$, $1$, $2$ and $3$, and that we would need to add the two extra solutions provided by block $4$ and its reversed one. Taking into account the subtlety mentioned, $2/3$ of times we would be one above the actual counting, so for those situations we only need a $+1$ ($-1$ for the surplus combination and $+2$ for blocks $4$ and the reversed one), whereas for the case when we are two below ($1/3$), we would need $+4$. The weighted number then is $+2$ to the total number obtained in (\ref{eq:Nclosedtotalcase12}) which leads to the final result in (\ref{eq:closeladsolutioncase1}).

\section{Circular ladder when $x/2>y$}\label{sec:appendix2}
In this case, the total number of allowed combinations for the circular ladder are:
\begin{equation}\label{eq:Nclosedtotal}
\begin{aligned}
    N_{circ}(k) = & N_{0..0}(k)+N_{0..3}(k)+N_{0..5}(k)+ N_{1..0}(k)+ \\ 
    & N_{1..3}(k)+N_{1..5}(k)+N_{2..1}(k)+N_{2..2}(k)+\\
    & N_{2..4}(k)+N_{3..1}(k)+N_{3..2}(k)+N_{3..4}(k)+\\
    & N_{4..0}(k)+N_{4..5}(k)+N_{5..2}(k)+N_{5..4}(k)
\end{aligned}
\end{equation}

We can rewrite each term as a function of $N_{0..0}(k)$ and $N_{0..2}(k)$ by using the relations in (42)-(45):
\begin{align}
    & N_{0..3}(k)=N_{0..2}(k) \\
    & N_{0..5}(k)=N_{0..2}(k)-N_{0..0}(k-1) \\
    & N_{1..0}(k)=N_{0..2}(k)-N_{0..0}(k-1) \\
    & N_{1..3}(k)=N_{0..0}(k)-N_{0..2}(k-1) \\
    & N_{1..5}(k)=N_{0..0}(k)+N_{0..0}(k-2)-2N_{0..2}(k-1) \\
    & N_{2..1}(k)=N_{0..2}(k) \\
    & N_{2..2}(k)=N_{0..0}(k) \\
    & N_{2..4}(k)=N_{0..2}(k)-N_{0..0}(k-1) \\
    & N_{3..1}(k)=N_{0..0}(k)-N_{0..2}(k-1) \\
    & N_{3..2}(k)=N_{0..2}(k)-N_{0..0}(k-1) \\
    & N_{3..4}(k)=N_{0..0}(k)+N_{0..0}(k-2)-2N_{0..2}(k-1) \\
    & N_{4..0}(k)=N_{0..2}(k) \\
    & N_{4..5}(k)=N_{0..0}(k)-N_{0..2}(k-1) \\
    & N_{5..2}(k)=N_{0..2}(k) \\
    & N_{5..4}(k)=N_{0..0}(k)-N_{0..2}(k-1) 
\end{align}\label{eq:long1}

Adding up all of those, (\ref{eq:Nclosedtotal}) becomes:
\begin{equation}\label{eq:closed0} 
\begin{aligned}
N_{circ}(k) = & 8\left[N_{0..0}(k)+N_{0..2}(k)\right] 
  -4\left[N_{0..0}(k-1)+N_{0..2}(k-1)\right] \\
 & +2\left[N_{0..0}(k-2)-2N_{0..2}(k-1)\right] 
\end{aligned}
\end{equation}

We also have:
\begin{equation}\label{eq:00and02}
    N_{0..0}(k)+N_{0..2}(k)=N_{0..}(k-1)=\frac{N(k-2)}{4}
\end{equation}
which when plugged into (\ref{eq:closed0}):
\begin{equation}\label{eq:Nalmostfinal}
\begin{aligned}
N_{circ}(k)=& 4N(k-2)-2N(k-3)+ \\
& +2[N_{0..0}(k-2)-2N_{0..2}(k-1)]
\end{aligned}
\end{equation}

And finally using $N_{0..0}(k)=N_{0..2}(k)=N(k-2)/4$ \footnote{As it was argued in section \ref{sec:appendix1} of the appendix, using this relation to simplify the third term in (\ref{eq:Nalmostfinal}) leads to the situation in which $2/3$ of the times we would be one above the real counting, and $1/3$ we would be two below, which translates into a weighted $0$ to fix the counting asymmetry.}, (\ref{eq:Nalmostfinal}) becomes:
\begin{equation}\label{eq:closed1}
N_{circ}(k)=4N(k-2)-3N(k-3)+\frac{N(k-4)}{2}
\end{equation}

From the recurrence relation in (\ref{eq:recurelation}) we can write $N(k-2)=3N(k-3)-N(k-4)$, then $3N(k-3)=N(k-2)+N(k-4)$, which when plugged into (\ref{eq:closed1}) finally gives us the result in (\ref{eq:closeladsolution}).

\newpage

\bibliography{biblioGT}                
\nocite{*}

\end{document}